%% file: ArxivDec2024.tex
\documentclass{amsart}

\usepackage[utf8]{inputenc}

\input{PreambleChapterProofs}

\author{Viktoriya Ozornova}
\address{Max Planck Institute for Mathematics, Bonn, Germany}
\email{viktoriya.ozornova@mpim-bonn.mpg.de}

\author{Martina Rovelli}
\address{Department of Mathematics and Statistics,
University of Massachusetts Amherst, 
Amherst,
USA
}
\email{mrovelli@umass.edu} 

\title{$(\infty,n)$-categories in context}

\begin{document}

\maketitle

\begin{abstract}
This note is a contribution written for the second volume of the \emph{Encyclopedia of mathematical physics}.
We give an informal introduction to the notions of an $(\infty,n)$-category and $(\infty,n)$-functor, discussing some of the different models that implement them. We also discuss the notions of a symmetric monoidal $(\infty,n)$-category and symmetric monoidal $(\infty,n)$-functor, recalling some important results whose statements employ the language of $(\infty,n)$-categories.
\end{abstract}

\section*{Introduction}

In mathematical physics, the language of higher categories is useful to propose formalisms that approximate physical theories and behaviors, such as topological quantum field theories, factorization algebras, $\mathbb{E}_n$-algebras, and higher gauge theory.

The presence of higher structures typically manifests itself whenever one tries to axiomatize the structure possessed by the class of certain mathematical objects and the expected axioms do not hold on the nose, but rather up to a ``higher isomorphism'' of some kind. One typical example that demonstrates this scenario is what happens when one tries to describe the mathematical aspects of the class of cobordisms: they can be naturally composed via appropriate gluing to produce new cobordisms, though the resulting composition law is only associative up to a boundary-preserving diffeomorphism.

We give a user-friendly overview of the notion of a (symmetric monoidal) $(\infty,n)$-category and $(\infty,n)$-functor, alongside examples that seem most relevant to situations from mathematical physics; we point to external references for a precise mathematical treatment when details are omitted.

In the first three sections we give an intuitive description of the notions of an $(\infty,n)$-category, $(\infty,n)$-functor and $(\infty,n)$-equivalence between $(\infty,n)$-categories. In the fourth section we discuss more in detail some of the mathematical models that formalize these concepts, describing some advantages and disadvantages of each. In the last three sections we discuss how the ideas of an $(\infty,n)$-category, an $(\infty,n)$-functor and an $(\infty,n)$-equivalence between $(\infty,n)$-categories interact with symmetric monoidal structures, and recall situations where the formalism is used to properly state theorems.

\addtocontents{toc}{\protect\setcounter{tocdepth}{1}}
\subsection*{Acknowledgements}
We are thankful to Lennart Meier for insightful discussions and to Catherine Meusburger and  Luuk Stehouwer for valuable feedback.
The second author is grateful for support from the National Science Foundation under Grant No. DMS-2203915.
\addtocontents{toc}{\protect\setcounter{tocdepth}{2}}
\section{$(\infty,n)$-Categories}
\label{ncat}

In this section we introduce the reader to the notion of an $(\infty,n)$-category from a heuristic viewpoint, and discuss a list of examples. Some more detailed proposed mathematical treatments of the notion of an $(\infty,n)$-category are postponed until the fourth section.

\subsection{Schematic definition}

We start by describing the type of structure expected from an $(\infty,n)$-category, and afterwards build towards the extra (structural) axioms that need to be added. This will consist of a type of categorical structure that has objects and morphisms in each positive dimension, with a weak composition of morphisms along lower dimensional morphisms, satisfying properties of weak associativity, weak unitality, and weak invertibility in dimensions higher than $n$.

\subsubsection*{The heuristic structure}
\label{basicstructure}
The structure of an \emph{$(\infty,n)$-category} consists of, for $k\geq0$, a set of \emph{$k$-morphisms} $k\operatorname{Mor}\mathscr{C}$, together with, for all $p\leq k$,
\begin{itemize}[leftmargin=*]
    \item \emph{source} and \emph{target} operators
\[s_p^k,t_p^k\colon k\operatorname{Mor}\mathscr{C}\to p\operatorname{Mor}\mathscr{C},\]
\item an \emph{identity} operator
\[\operatorname{id}_k^p\colon p\operatorname{Mor}\mathscr{C}\to k\operatorname{Mor}\mathscr{C},\]
\item and a \emph{composition} operator
\[\circ_p^k\colon k\operatorname{Mor}\mathscr{C}\times_{p\operatorname{Mor}\mathscr{C}}k\operatorname{Mor}\mathscr{C}\to k\operatorname{Mor}\mathscr{C},\]
 defined on all pairs of $k$-morphisms $(f,g)$ for which $t^k_p(f)=s^k_p(g)$.
\end{itemize}

By convention, a $0$-morphism is also referred to as an \emph{object}, a $1$-morphism is referred to as a \emph{morphism}, and a $k$-morphism for $k>1$ is referred to as a \emph{higher morphism}. The \emph{dimension} of a $k$-morphism is $k$.

For ease of exposition, superscripts and subscripts are often dropped, so one could write $s(f)$ and $t(f)$ in place of $s_{k-1}^k(f)$ and $t_{k-1}^k(f)$, or $\operatorname{id}_f$ for $\operatorname{id}_{k+1}^k(f)$, or $g\circ_{p}f$ for $g\circ^k_{p}f$, or $g\circ f$ for $g\circ_{k-1}^kf$.
If $s(f)=a$ and $t(f)=b$ for a $k$-morphism $f$ of $\mathscr{C}$, we depict $f\colon a\to b$.

\subsubsection*{Replacing axioms with higher equivalences}

The structure of an $(\infty,n)$-category $\mathscr{C}$ needs to be complemented in order to encode appropriately the facts that, in some sense:
\begin{itemize}[leftmargin=*]
    \item the various compositions are associative;
    \item the various compositions are unital, with identity morphisms of lower-dimensional morphisms being the neutral elements;
    \item the various compositions satisfy certain interchange laws;
    \item every morphism in dimension higher than $n$ is invertible with respect to composition.
\end{itemize}

Normally we'd give a list of axioms for this purpose, but this would produce a too restrictive notion of an $(\infty,n)$-category, failing to recover examples of interest. Given an axiom of interest, rather than imposing it in the form of an equality between certain $k$-morphisms in $\mathscr{C}$, we shall require the existence of further structure -- a ``$(k+1)$-equivalence'' in $\mathscr{C}$ between such $k$-morphisms. In this case we think of the $(k+1)$-equivalence as a \emph{witness} -- witnessing
for the relation that enhances the naive axiom appropriately. We may think of these as \emph{structural axioms} for $\mathscr{C}$.

We define by a downward induction what it means for a $k$-morphism $f\colon a\to b$ to be a $k$-equivalence in $\mathscr{C}$, in which case we'll denote $f\colon a\simeq b$.
We write $a\simeq b$ when there exists a $k$-equivalence $f\colon a\simeq b$.

\begin{itemize}[leftmargin=*]
    \item For $k>n$, every $k$-morphism of $\mathscr{C}$ is to be considered a \emph{$k$-equivalence}.
    \item For $k\leq n$ a $k$-morphism $f\colon a\to b$ of $\mathscr{C}$ is to be considered a \emph{$k$-equivalence} when there exists a $k$-morphism $g\colon b\to a$ and $(k+1)$-equivalences of $\mathscr{C}$
\[\operatorname{id}_b\simeq f\circ g\text{ and }\operatorname{id}_a\simeq g\circ f.\]
\end{itemize}

\subsubsection*{The first layer of structural axioms}
\label{The first layer of structural axioms}

The most superficial layer of the structural axioms that we shall impose on the structure of $(\infty,n)$-category $\mathscr{C}$ resembles the traditional list of axioms that one would write for a \emph{strict} $n$-category, upon replacing each equality sign with the existence of a higher equivalence in $\mathscr{C}$.
\begin{itemize}[leftmargin=*]
    \item \emph{Associativity:} Given $k$-morphisms
$f\colon a\to b$, $g\colon b\to c$ and $h\colon c\to d$
    there is a $(k+1)$-equivalence
\[h\circ (g\circ f)\simeq(h\circ g)\circ f.\]
\item \emph{Unitality:} Given a $k$-morphism $f\colon a\to b$, there are $(k+1)$-equivalences
\[f\circ\operatorname{id}\simeq f\quad\text{ and }\quad\operatorname{id}\circ f\simeq f\]
\item \emph{Interchange:} Given $k$-morphisms
$f_i\colon a\to b$ and $g_i\colon b\to c$
for $i=1,2,3$, and given $(k+1)$-morphisms
$\alpha_i\colon f_i\to f_{i+1}$ and $\beta_i\colon g_i\to g_{i+1}$
for $i=1,2$, there is a $(k+2)$-equivalence
\[(\alpha_2\circ_{k-1}\alpha_1)\circ_{k-2}(\beta_2\circ_{k-1}\beta_1)\simeq(\beta_2\circ_{k-2}\alpha_2)\circ_{k-1}(\beta_1\circ_{k-2}\alpha_1).\]
\item \emph{Invertibility:} Given a $k$-morphism $f\colon a\to b$ for $k>n$ there exists a $k$-morphism $g\colon b\to a$ and $(k+1)$-equivalences
\[\operatorname{id}_a\simeq g\circ f\quad\text{ and }\quad f\circ g\simeq \operatorname{id}_b.\]
\end{itemize}

\subsubsection*{A taste for higher coherence}
\label{A taste for higher coherence}

The list of structural axioms that we proposed for the $(\infty,n)$-category $\mathscr{C}$ is not enough. Indeed, replacing traditional equalities with the existence of witnessing higher equivalences forces as well the existence of infinitely many new expressions that one could build out of such higher equivalences, and it is necessary to impose further coherence between them.

To showcase this principle, let's just have a closer look at the associativity requirement; the others would present similar issues. To this end, let's denote
\[\alpha_{f,g,h}\colon h\circ (g\circ f)\simeq(h\circ g)\circ f\]
any $2$-equivalence witnessing the associativity relation for a generic triple of $1$-morphisms $f$, $g$, and $h$. Then, given four composable $1$-morphisms
\[f\colon a_0\to a_1\quad g\colon a_1\to a_2\quad h\colon a_2\to a_3\quad \text{ and } \quad\ell\colon a_3\to a_4,\]
there are at least two ways to produce a
$2$-morphism between the $1$-morphisms $(\ell\circ (h\circ(g\circ f)))\colon a_0\to a_4$ and $(((\ell\circ h)\circ g)\circ f)\colon a_0\to a_4$:
the $2$-morphism
\[(\alpha_{g,h,\ell}\circ_0\operatorname{id}_{f})\circ \alpha _{f, h \circ g, \ell}\circ (\operatorname{id}_{\ell}\circ_0\alpha_{f, g, h})\colon\ell\circ (h\circ(g\circ f))\to((\ell\circ h)\circ g)\circ f\]
and the $2$-morphism
\[\alpha_{f, g, \ell\circ h} \circ \alpha_{g\circ f, h, \ell}\colon \ell\circ (h\circ(g\circ f))\to((\ell\circ h)\circ g)\circ f.\]
A necessary coherence requirement that needs to be imposed is then the existence of a $3$-equivalence
\[\xi_{f,g,h,\ell}\colon(\alpha_{g,h,\ell}\circ_0\operatorname{id}_{f})\circ \alpha _{f, h \circ g, \ell}\circ (\operatorname{id}_{\ell}\circ_0\alpha_{f, g, h})\simeq\alpha_{f, g, \ell\circ h} \circ \alpha_{g\circ f, h, \ell}.\]
As one might expect, the $3$-equivalence $\xi_{f,g,h,\ell}$ would then also in turn need to be subject to coherence conditions, showcasing how the infinite amount of required coherence datum arises, even just looking at associativity in the lowest dimension.
\[
\begin{tikzcd}[column sep=1.1cm, /tikz/column 2/.style={column sep=1.5cm}]
& \ell\circ ((h\circ g) \circ f)\arrow[r, "\alpha _{f, h \circ g, \ell}", ""{name=T,inner sep=2pt, swap}]  & (\ell\circ (h\circ g)) \circ f\arrow[dd, "\alpha_{g,h,\ell}\circ_0\operatorname{id}_{f}"]&\\
\ell\circ (h\circ(g\circ f)) \arrow[ur, "\operatorname{id}_{\ell}\circ_0\alpha_{f, g, h}"] \arrow[dr, "\alpha_{g\circ f, h, \ell}" swap]&&&\\
&(\ell \circ h) \circ (g \circ f) \arrow[r, "\alpha_{f, g, \ell\circ h}" swap]&((\ell\circ h)\circ g)\circ f&
\arrow[from=T, to=3-2, shorten >= 0.6cm, shorten <= 0.3cm, "\xi_{f,g,h,\ell}"]
\end{tikzcd}
\]

\subsubsection*{About composition being structural}

Beside higher coherence, a somewhat orthogonal issue is that there is a sense in which it is unnatural and not completely appropriate to assume the composition operators (as well as the identity, source and target operators) to be well-defined functions included in the structure of $\mathscr{C}$. Instead, it would be more appropriate to require that, given any composable configuration of morphisms in $\mathscr{C}$, there must exist a composite alongside a witness for this composite.
As a first order approximation, and for ease of exposition, we maintain the existence of the composition functions as part of the structure of an $(\infty,n)$-category.
It is likely that, anyway, one could with loss of generality assume the existence of structural composition maps. For instance, the arguments of \cite{NikolausAlgebraicModels} should provide useful tools to formalize this in the model of $(\infty,n)$-categories given by $n$-quasi-categories \cite{ara}.

\subsubsection*{Different models for an $(\infty,n)$-category}

\label{Different models}

In total, the structure of an $(\infty,n)$-category consists of the basic heuristic structure that we just discussed, as well as lots of higher coherence data which is essentially impossible to express as a list.
Given the infinite amount of coherence that one needs to keep track of, it should already be clear at this point that it would not be manageable to give a precise definition of an $(\infty,n)$-category based on the usual list of structure maps and axioms.

Several mathematical frameworks have been introduced to formalize precisely the definition of an $(\infty,n)$-category, giving meaning to the notions of $k$-morphisms and $k$-equivalences, while packaging the appropriate coherence data in different ways. These include:
\begin{itemize}[leftmargin=*]
    \item \emph{Models based on various flavors of marked simplicial sets:} These include Kan complexes \cite[\textsection II.3]{QuillenHA} for $n=0$, quasi-categories \cite{JoyalVolumeII} and naturally marked quasi-categories
    \cite[\textsection 3.1]{htt} for $n=1$, $\infty$-bicategories \cite{LurieGoodwillie,GHL} for $n=2$, and $n$-complicial sets \cite{VerityComplicialI,or}.
    \item \emph{Various flavors of (marked) cubical sets:} These include $n$-comical sets \cite{DKLS,CKM,DKM}.
    \item \emph{Models based on multisimplicial spaces:} These include complete Segal spaces \cite{rezk} for $n=1$, $n$-fold complete Segal spaces \cite{BarwickThesis} and categories weakly enriched in $(\infty,n-1)$-categories \cite{br1}.
    \item \emph{Models based on globular presheaves:} These include complete Segal $\Theta_n$-spaces \cite{rezkTheta}, $n$-quasi-categories \cite{ara}, and the colossal model from \cite{BarwickSchommerPries}.
    \item \emph{Models based on categories strictly enriched over models $(\infty,n-1)$-categories:} These include the enriched models from \cite{br1}.
\end{itemize}
We discuss some of them -- $n$-complicial sets, $n$-fold complete Segal spaces, and categories strictly enriched over $(\infty,n-1)$-categories -- later, and refer the reader to e.g.~\cite{BergnerModelsn,bergnermodels} for a survey on some of the most popular different approaches.
By now it is known that all the above approaches are equivalent in an appropriate sense, and we refer the reader to
\cite{ara,br1,br2,DKM,LoubatonEqui,LurieGoodwillie}
for sources that discuss some of the comparisons.

\subsection{Examples}

\label{ExamplesCategories}

We discuss several examples of $(\infty,n)$-categories that appear naturally for various values of $n$. For many of the examples multiple (and substantially different!) variants have been considered in the literature, and we invite the reader to keep this in mind to avoid confusion.

\subsubsection*{Structural examples of $(\infty,n)$-categories}

We collect here examples of $(\infty,n)$-ca\-te\-go\-ries that arise from other known variants of (higher) categories.

\begin{ex}
An $(\infty,n)$-category can be naturally regarded as an $(\infty,n+1)$-category. Given an $(\infty,n+1)$-category $\mathscr{D}$, there are two universal ways of producing an $(
\infty,n)$-category: the core $\mathrm{core}_n\mathscr{D}$, in which one forgets the non-invertible $(n+1)$-morphisms, and the intelligent truncation $\tau\mathscr{D}$, which freely inverts all $(n+1)$-morphisms. It is recalled e.g.~in \cite[\textsection2.1.2]{ORsurvey} how to formalize this construction in the model of $n$-complicial sets, based on constructions from \cite{VerityComplicialI}.
\end{ex}

\begin{ex}
Given an $(\infty,n)$-category $\mathscr{C}$ and two objects $a$ and $b$ in $\mathscr{C}$, there is an $(\infty,n-1)$-category $\mathscr{C}(a,b)$, called the \emph{hom-$(\infty,n-1)$-category} in which the objects are the morphisms $f\colon a\to b$ in $\mathscr{C}$, and a $k$-morphism in $\mathscr{C}(a,b)$ corresponds to a $(k+1)$-morphism in $\mathscr{C}$ for which the $0$-dimensional source and target are $a$ and $b$, respectively. It is recalled in \cite[\textsection2.1.2]{ORsurvey} (based on \cite{VerityComplicialI}) how to formalize this construction for the model of $n$-complicial sets, in \cite{rezkTheta} for the model of complete Segal $\Theta_n$-spaces, and in \cite[\textsection A.3]{htt} and \cite[\textsection3.10]{br1} for the model of categories enriched over $(\infty,n-1)$-categories.
\end{ex}

\begin{ex}
\label{smallnerves}
A category can be regarded naturally as an $(\infty,1)$-category. A $2$-category or more generally a bicategory in the sense of \cite[\textsection 1]{Benabou} can be regarded as an $(\infty,2)$-category. This is formalized in \cite{CampbellHoCoherent} in the model of $2$-quasi-categories.
\end{ex}

\begin{ex}
A topological space (which we refer to simply as a space) can be naturally regarded as an $(\infty,0)$-category, as we'll describe in \cref{ExampleSpace}. In fact, there is a sense in which all $(\infty,0)$-categories are spaces. This is formalized e.g.~in \cite[\textsection II.3]{QuillenHA} in the model of Kan complexes.
\end{ex}

\subsubsection*{$(\infty,n)$-categories of structured sets}
We give a family of examples of $(\infty,n)$-categories in which the objects are structured sets and the morphisms are appropriate structure-preserving functions. In these cases, composition amounts to usual composition of functions and identity amounts to the usual identity function.

\begin{ex}[Vector spaces \& Hilbert spaces]
Given a field $\mathbb{K}$, there is an $(\infty,1)$-category -- in fact an ordinary category -- $\mathscr{V}ect_{\mathbb{K}}$ (resp.~$\mathscr{H}ilb_{\mathbb{K}}$, $\mathscr{A}lg_{\mathbb{K}}$) in which the objects are $\mathbb{K}$-vectors spaces (resp.~$\mathbb{K}$-Hilbert spaces, $\mathbb{K}$-algebras), the morphisms are $\mathbb{K}$-linear maps (resp.~$\mathbb{K}$-bounded operators, $\mathbb{K}$-algebra maps), and there are no non-trivial higher morphisms.
\end{ex}

\begin{ex}[Spaces]
There is an $(\infty,1)$-category -- which is \emph{not} an ordinary category -- $\mathscr{S}pc$ in which the objects are spaces, the morphisms are continuous functions, and the higher morphisms are homotopies, homotopies between homotopies, and so on. This can be formalized as the underlying $(\infty,1)$-category of the model category from \cite[\textsection II.3]{QuillenHA}.
\end{ex}

\begin{ex}[Chain complexes]
Given a field $\mathbb{K}$, there is an $(\infty,1)$-category -- which is \emph{not} an ordinary category -- $\mathscr{C}h_{\mathbb{K}}$ in which the objects are chain complexes over $\mathbb{K}$, the morphisms are chain maps, the $2$-morphisms are chain homotopies, and the higher morphisms are higher chain homotopies between chain homotopies in an appropriate sense.
This can be formalized as the underlying $(\infty,1)$-category of the model category from \cite[\textsection 2.3]{hovey}.
\end{ex}

\begin{ex}[Linear categories]
Given a field $\mathbb{K}$, there is an $(\infty,2)$-category -- in fact a $2$-category --
$\mathscr{C}at_{\mathbb{K}}^L$,
where the objects are finitely cocomplete linear categories, and the $1$-morphisms are the right exact functors. As $2$-morphisms one could take all the $\mathbb{K}$-linear natural transformations between right exact
functors (cf.~\cite[(3.1)]{TelemanLectures}) or select only the invertible ones (cf.~\cite[\textsection 3]{BenZviBrochierJordanIntegrating}), and there are no non-trivial higher morphisms.
\end{ex}

\subsubsection*{$(\infty,n)$-categories with unusual morphisms}
\label{UnusualMorphisms}

It is not always the case that a morphism in an $(\infty,n)$-category must be a function preserving some structure, and here are some (less intuitive!) examples that showcase this.

\begin{ex}[Points in a space]
\label{ExampleSpace}
Given a space $X$, there is an $(\infty,0)$-category $\mathscr{X}$ in which the objects are the points of $X$, the morphisms $\alpha\colon x\to y$ correspond to paths from $x$ to $y$
\[\alpha\colon x\to y\text{ in }\mathscr{X}\quad\leftrightsquigarrow\quad \alpha\colon[0,1]\to X\text{ in }\mathscr{S}pc\text{ with }\left\{
\begin{array}{cc}
    \alpha(0)=x &  \\
    \alpha(1)=y & 
\end{array}\right.\]
and the higher morphisms are given by homotopies of paths in $X$, homotopies of homotopies, and so on.
In this context, composition is given by appropriate concatenation of paths or homotopies in $X$, and identity is given by constant paths in $X$. The $(\infty,0)$-category $\mathscr{X}$ is essentially formalized in the form of a Kan complex by the singular simplicial set construction of $X$.
\end{ex}

\begin{ex}[Spaces and Spans]
There is an $(\infty,n)$-category
$\mathscr{S}pan_n$ in which the objects are spaces, sets, or groupoids,
possibly satisfying further finiteness conditions, the $1$-morphisms $\phi\colon X\to Y$ correspond to spans
\[\phi\colon X\to Y\text{ in }\mathscr{S}pan_{n}\quad\leftrightsquigarrow\quad\left[X\xleftarrow{\phi_1} S\xrightarrow{\phi_2} Y\right]\text{ in }\mathscr{S}pc.\]
There are several choices for the $2$-morphisms and higher morphisms in general. One could take as $k$-morphisms in $\mathscr{S}pan_n$ for $k\leq n$ either $k$-fold iterated spans or maps of spans, and then as $k$-morphisms for $k>n$ weak equivalences or isomorphisms of spans. The former is the viewpoint taken e.g.~in \cite[\textsection3]{barwickq} for $n=1$,
and \cite[\textsection2.6]{Benabou} and \cite[\textsection10]{DK} for $n=2$; see also \cite{CalaqueHaugsengScheimbauer}. The latter is the viewpoint taken e.g.~in \cite[\textsection3.2]{luriecobordism} and  \cite[\textsection5]{HaugsengSpan}.
In this context, composition of spans is given by (homotopy-)pulling back appropriately.
\end{ex}

\begin{ex}[Algebras and Bimodules]
Given a field $\mathbb{K}$,
there is an $(\infty,2)$-category -- in fact a bicategory -- $\mathscr{M}or_{\mathbb{K}}$ -- sometimes referred to as \emph{Morita bicategory} -- in which the objects are $\mathbb{K}$-algebras, the $1$-morphisms from $R$ to $S$ are $(R,S)$-bimodules
\[M\colon R\to S\text{ in }\mathscr{M}or_{\mathbb{K}}\quad\leftrightsquigarrow\quad {}_RM_S=\left[R \curvearrowright M \curvearrowleft S\right]\text{ bimodule},\]
the $2$-morphisms are given by bimodule morphisms and there are no interesting higher morphisms. In this context, composition along a $0$-morphism is given by tensoring over the appropriate algebra, and composition along a $1$-morphism is ordinary composition of bimodule maps. Identity of objects is given by regarding a $\mathbb{K}$-algebra as a bimodule over itself, and identity of a $1$-morphism is given by the usual identity morphism of a bimodule. The $(\infty,2)$-category $\mathscr{M}or_{\mathbb{K}}$ is formalized in the form a bicategory e.g.~in \cite[\textsection 4.6]{SPthesis}.
\end{ex}

\subsubsection*{$(\infty,n)$-categories from a group}

We mention two types of higher categories arising naturally from a (finite, discrete, topological or Lie) group $G$. Various groups (e.g.\ symmetry groups and unitary groups) play a role in formalizing many physical situations, and a starting point of reading for their occurrence e.g.\ in the contexts involving TQFTs could be \cite{Snowmass2022}

\begin{ex}
Given a (finite, discrete, topological or Lie) group $G$, there is an $(\infty,0)$-category $\mathscr{B}G$ in which there is exactly one object called $\ast$, the morphisms $g\colon *\to *$ correspond to the elements $g$ of $G$
\[g\colon*\to*\text{ in }\mathscr{B}G\quad\leftrightsquigarrow\quad g\in G,\]
and higher morphisms are given by paths, homotopies, homotopies between homotopies and so on in $G$. In this context, composition and identity are given by the operation and neutral element of $G$.
The $(\infty,0)$-category $\mathscr{B}G$ can be formalized as a topological category with one object or as a Kan complex.
\end{ex}

\begin{ex}[Principal bundles]
Given a space $X$ and a finite\footnote{In the case of $G$ being infinite the same construction makes sense but it would not have many of the desirable properties (for instance it would not satisfy the analog statement to \cref{ClassificationGBundles}). In that case, one should instead consider an appropriate $(\infty,0)$-category of ``weak principal bundles'', as in \cite{NikolausSchreiberStevensonPresentations}.} group $G$, there is an 
$(\infty,0)$-category $\mathscr{B}un^G_X$ in which the objects are the principal $G$-bundles over $X$, the morphisms are the $G$-equivariant maps of bundles over $X$, and there are no non-trivial higher morphisms.
\end{ex}

\subsubsection*{$(\infty,n)$-categories of manifolds}

\label{ExamplesManifolds}

We collect a list of examples of $(\infty,n)$-categories whose objects are manifolds of some kind. 

We will discuss in the following sections their role,
respectively, in notions such as \emph{factorization algebra} \cite[\textsection 6.1]{CostelloGwilliam1}, \cite[\textsection5.5.2]{LurieHA}, \emph{$\mathbb{E}_n$-algebra} \cite[\textsection 2.7]{AFprimer}, \cite[\textsection 5.1.1]{LurieHA},
\emph{factorization homology} \cite{AF1,AFprimer} and \emph{topological chiral homology} \cite[\textsection5.5.2]{LurieHA}
, and (extended) \emph{topological quantum field theory} \cite{AtiyahTQFT,luriecobordism}.

\begin{ex}[Open subsets]
Given $M$ a manifold
(e.g.~ $M=\mathbb{R}^n$), there is an $(\infty,1)$-category -- in fact an ordinary category -- $\mathscr{O}p_{M}$, in which the objects are open subsets of $M$, the $1$-morphisms are inclusions, and there are no non-trivial higher morphisms.
\end{ex}

\begin{ex}[Manifolds, Disks]
Given $n\geq0$
, there is an $(\infty,1)$-category $\mathscr{M}fld^{\mathrm{fr}}_{n}$ where, roughly speaking, the objects are framed $n$-dimensional
smooth manifolds that satisfy an appropriate finiteness condition\footnote{The finiteness condition requires the existence of a finite good cover, which is to say a finite open cover by Euclidean spaces with the property that each non-empty intersection of terms in the cover is itself homeomorphic to a Euclidean space.} and the $1$-morphisms are 
embeddings 
which preserve framings in a suitable homotopical sense.
The higher morphisms essentially encode isotopies, isotopies of isotopies, and so on.
The $(\infty,1)$-category $\mathscr{M}fld^{\mathrm{fr}}_{n}$ is formalized precisely as a topological category in \cite[\textsection 2.4-2.5]{AFprimer}. The $(\infty,1)$-category $\mathscr{D}isk^{\mathrm{fr}}_{n}$ is defined similarly, with the objects being disjoint unions of finitely many copies of a framed
$n$-dimensional euclidian space $\mathbb{R}^n$, and is formalized as a topological category in \cite[\textsection 2.7]{AFprimer}.
\end{ex}

\begin{ex}[Cobordisms]
Given $n>0$,
there is an $(\infty,1)$-category $\mathscr{C}ob_{n-1,n}$ -- in fact a category -- where, roughly speaking, the objects are
closed $(n-1)$-manifolds,
the morphisms are $n$-dimensional manifolds
seen as a cobordism between their incoming and their outgoing boundary components
\[\phi\colon U\to V\text{ in }\mathscr{C}ob_{n-1,n}\quad\leftrightsquigarrow\quad\left[U\stackrel{\phi_1}\hookrightarrow \Sigma\stackrel{\phi_2}\hookleftarrow V\right]\quad\text{s.t.}\quad(\phi_1,\phi_2)\colon U\amalg V\cong\partial\Sigma\]
up to boundary preserving diffeomorphisms, and there are no non-trivial higher morphisms.
In this context composition is given by gluing cobordisms and composing diffeomorphism classes, and identity is given by taking the cylinder cobordism. A version of the $(\infty,1)$-category $\mathscr{C}ob_{n-1,n}$ is formalized as a category e.g.~in
\cite[\textsection2.1]{GMTWCob}.
 Variants of the category $\mathscr{C}ob_{n-1,n}$ are defined similarly but take into account manifolds and cobordism with extra geometric structure (such as orientation, framing, conformal structure, etc.), which is recorded in the notation (e.g.,~$\mathscr{C}ob^{\mathrm{or}}_{n-1,n}$, $\mathscr{C}ob^{\mathrm{fr}}_{n-1,n}$, $\mathscr{C}ob^{\mathrm{conf}}_{n-1,n}$).
\end{ex}

\begin{ex}[Cobordisms, extended]
Given $n\geq0$, there is an $(\infty,n)$-category $\mathscr{C}ob_{0,n,\infty}$. Heuristically, the objects are
closed $0$-manifolds, for $0<k<n+1$ the $k$-morphisms are $k$-dimensional manifolds with corners seen as a cobordism between their incoming and their outgoing boundary components,
$(n+1)$-morphisms are boundary preserving diffeomorphisms, and higher morphisms are isotopies of such. Variants of the category $\mathscr{C}ob_{0,n,\infty}$ are defined similarly but take into account manifolds and cobordisms with extra geometric structure, and\kern-0.001cm /\kern-0.001cm or focus on different dimension ranges for the involved cobordisms. For instance, a version of the $(\infty,n)$-category $\mathscr{C}ob^{\mathrm{fr}}_{0,n,\infty}$ is formalized precisely in \cite[\textsection 5]{CalaqueScheimbauer} in the form of an $n$-fold complete Segal space; see also \cite[\textsection 2.2]{luriecobordism}. A version of a bicategory $\mathscr{C}ob_{0,2}^{\mathrm{or}}$ is constructed in \cite[\textsection4]{SPthesis}, and a version of the $(\infty,2)$-category $\mathscr{C}ob^{\mathrm{or}}_{1,3}$ is described in \cite{BDSPV}.
\end{ex}

\section{$(\infty,n)$-Functors}
\label{nfun}

In this section we introduce the reader to the notion of an $(\infty,n)$-functor between $(\infty,n)$-categories from a heuristic viewpoint, and discuss a list of examples. Some of the mathematical treatments of the notion of an $(\infty,n)$-functor
are postponed until later.

\subsection{Schematic definition}

We start by describing the type of structure expected from an $(\infty,n)$-functor, and afterwards build towards the extra (structural) axioms that need to be added.

\subsubsection*{The heuristic structure}
Given $(\infty,n)$-categories $\mathscr{C}$ and $\mathscr{D}$, the structure of an \emph{$(\infty,n)$-functor} $F\colon\mathscr{C}\to\mathscr{D}$
consists of an operator
\[F_k\colon k\operatorname{Mor}\mathscr{C}\to k\operatorname{Mor}\mathscr{D}\]
for all $k\geq0$. Like for $(\infty,n)$-categories, we will often write $F$ in place of $F_k$.

We shall require the existence of further structure to witness an appropriate enhancement of the usual axioms we would require for a functor, asserting compatibility with the source and target operators, the composition operators, and the identity operators of $\mathscr{C}$ and of $\mathscr{D}$.

\subsubsection*{The first layer of structural axioms}

\label{The first layer of structural axioms 2}

The first layer of requirements for an $(\infty,n)$-functor $F\colon\mathscr{C}\to\mathscr{D}$ 
amounts to the following conditions for the structure maps of $\mathscr{C}$ and $\mathscr{D}$:
\begin{itemize}[leftmargin=*]
    \item \emph{Source and target:} Given a $k$-morphism $f$, there is an equality of $k$-morphisms
\[s(F(f))= F(s(f))\quad\text{ and }\quad t(F(f))= F(t(f));\]
\item \emph{Identity:} Given a $k$-morphism $f$
there is a $(k+2)$-equivalence
\[F(\operatorname{id}_f)\simeq \operatorname{id}_{Ff};\]
\item \emph{Composition:}
Given $k$-morphisms $f\colon a\to b$ and $g\colon b\to c$ there is a $(k+1)$-equivalence
\[F(g\circ f)\simeq F(g)\circ F(f).\]
\end{itemize}

\subsubsection*{A taste for higher coherence}

\label{A taste for higher coherence 2}
In a similar vein to what we discussed in the first section, the list of axioms that we proposed for the
$(\infty,n)$-functor $F\colon\mathscr{C}\to\mathscr{D}$ is not enough, and we describe one manifestation of the issue in this context.

Let's denote
\[\gamma_{f,g}\colon F(g\circ f)\simeq F(g)\circ F(f)\]
any $2$-equivalence witnessing the functoriality relation for a generic pair of $1$-morphisms $f$ and $g$ in $\mathscr{C}$. Let's denote by $\alpha^{\mathscr{C}}_{f,g,h}$ and $\alpha^{\mathscr{D}}_{Ff,Fg,Fg}$
the equivalences in $\mathscr{C}$ and $\mathscr{D}$ for three given composable $1$-morphisms
\[f\colon a_0\to a_1\quad g\colon a_1\to a_2\quad h\colon a_2\to a_3.\]
There are then at least two ways to produce a 
$2$-morphism in $\mathscr{D}$ between the $1$-morphisms $F(h\circ (g\circ f))\colon Fa_0\to Fa_3$ and $(Fh\circ Fg)\circ Ff\colon Fa_0\to Fa_3$: the $2$-morphism in $\mathscr{C}$
\[(\gamma_{g,h}\circ \operatorname{id}_{Ff})\circ\gamma_{f, h\circ g}\circ F(\alpha^{\mathscr{C}}_{f,g,h})\colon F(h\circ (g\circ f))
\to (Fh\circ Fg)\circ Ff\]
and the
$2$-morphism
\[\alpha^{\mathscr{D}}_{Ff, Fg, Fh}\circ (\operatorname{id}_{Fh}\circ \gamma_{f,g})\circ \gamma_{g\circ f, h}\colon F(h\circ (g\circ f))
\to(Fh\circ Fg)\circ Ff.\]
A further necessary coherence requirement is the existence of a $3$-equivalence in $\mathscr{D}$
\[\nu_{f,g,h}\colon(\gamma_{g,h}\circ \operatorname{id}_{Ff})\circ\gamma_{f, h\circ g}\circ F(\alpha^{\mathscr{C}}_{f,g,h})\simeq \alpha^{\mathscr{D}}_{Ff, Fg, Fh}\circ (\operatorname{id}_{Fh}\circ \gamma_{f,g})\circ\gamma_{g\circ f, h},\]
and there are of course (infinitely) many more.
\[
\begin{tikzcd}[column sep=0.7cm, /tikz/column 2/.style={column sep=1.5cm}]
& F((h\circ g)\circ f)\arrow[r, "\gamma_{f, h\circ g}", ""{name=T,inner sep=2pt, swap}]  & F(h\circ g) \circ Ff\arrow[dr, "\gamma_{g,h}\circ \operatorname{id}_{Ff}"]&\\
F(h\circ (g\circ f)) \arrow[ur, "F(\alpha^{\mathscr{C}}_{f,g,h})"] \arrow[dr, "\gamma_{g\circ f, h}" swap]&&& (Fh\circ Fg)\circ Ff\\
&Fh \circ F(g\circ f) \arrow[r, "\operatorname{id}_{Fh}\circ \gamma_{f,g}" swap, ""{name=B,inner sep=2pt} ]&Fh\circ (Fg \circ Ff) \arrow[ur, "\alpha^{\mathscr{D}}_{Ff, Fg, Fh}" swap]&
\arrow[from=T, to=B, shorten >= 0.6cm, shorten <= 0.3cm, "\nu_{f,g,h}"]
\end{tikzcd}
\]

Each of the models proposed to model the notion of an $(\infty,n)$-category must come with a corresponding notion of an $(\infty,n)$-functor that packages the appropriate amount of coherence. Various explicit implementations will be discussed later.

\subsection{Examples}

We discuss several examples of $(\infty,n)$-functors that appear naturally for various values of $n$.

\subsubsection*{Structural examples}

We collect here examples of $(\infty,n)$-functors that arise from other known variants of functors between (higher) categories.

\begin{ex}
Every $(\infty,n)$-functor can naturally be regarded as an $(\infty,n+1)$-functor.
\end{ex}

\begin{ex}
Every ordinary functor between ordinary categories defines an $(\infty,1)$-functor, and every normal pseudo-functor (as in \cite[\textsection4]{Benabou}, there referred to as \emph{unitary homomorphism})
is an $(\infty,2)$-functor (cf.~\cref{smallnerves}).
\end{ex}

\begin{ex}
Every map of spaces defines an $(\infty,0)$-functor, and there is a sense in which every $(\infty,0)$-functor is a map of spaces.
\end{ex}

\begin{ex}
Given an $(\infty,n)$-category $\mathscr{C}$, there is an \emph{identity $(\infty,n)$-functor} $\operatorname{id}_\mathscr{C}\colon\mathscr{C}\to\mathscr{C}$.
\end{ex}

\begin{ex}
Given $(\infty,n)$-functors $F\colon\mathscr{C}\to\mathscr{D}$ and $G\colon\mathscr{D}\to\mathscr{E}$, there is a \emph{composite $(\infty,n)$-functor} $G\circ F\colon\mathscr{C}\to\mathscr{E}$.
\end{ex}

\subsubsection*{Explicit examples}

We see some examples that involve some of the explicit examples from \cref{ExamplesCategories} which are often considered. We warn the reader that, while all the $(\infty,n)$-functors are surely formalizable in an appropriate context, one has to be very careful in choosing the correct incarnation and flavor of the $(\infty,n)$-categories involved, and -- even then -- it is not always possible to find a rigorous implementation in the existing literature.

\begin{ex}
Let $X$ be a topological space and $G$ a group.
An $(\infty,0)$-functor
\[f\colon\mathscr{X}\to \mathscr{B}G\]
has no interesting information
in dimension $0$, and
consists of a coherent way of assigning an element in $G$ to each path in $X$, so that homotopic paths get sent to the same element, and a concatenation of paths gets sent to the product of the corresponding elements in $G$.
We will see in \cref{ClassificationGBundles} how this is the same information as encoded 
into $\mathscr{B}un^G_X$ at least when $G$ is finite.
\end{ex}

\begin{ex}
\label{ExampleBunFunctor}
Let $G$ be a group. For appropriate choices of the $(\infty,2)$-categories $\mathscr{C}ob_{0,2,\infty}^{\mathrm{or}}$ and $\mathscr{S}pan_2$, there is an $(\infty,2)$-functor
\[\mathscr{B}un^G_{(-)} \colon\mathscr{C}ob_{0,2,\infty}^{\mathrm{or}}\to\mathscr{S}pan_2,\]
which is essentially given on objects by
\[U\mapsto \mathscr{B}un^G_U,\]
on $k$-morphisms for $k=1,2$ by
\[[U\stackrel{\phi}\hookrightarrow\Sigma\stackrel{\phi'}\hookleftarrow  U']\mapsto\left[\mathscr{B}un^G_U\stackrel{\phi^*}\longleftarrow \mathscr{B}un^G_\Sigma\stackrel{(\phi')^*}\longrightarrow \mathscr{B}un^G_{U'}\right].\]
A version of this $(\infty,2)$-functor is discussed in \cite[\textsection3]{TelemanLectures} and as the case $n=2$ of \cite[(3.3)]{FHLT}.
\end{ex}

\begin{ex}
\label{ExampleCategorification}
Let $\mathbb{K}$ be a field. There is an $(\infty,2)$-functor
\[\mathscr{M}or_\mathbb{K}\to\mathscr{C}at^L_{\mathbb{K}},\]
which is given on objects by
\[A\mapsto \mathcal{M}od_{A}\]
and on morphisms by
\[{}_AM_B\mapsto \left[-\otimes_AM_B\colon\mathcal{M} od_A\to\mathcal{M} od_B\right]\]
and on $2$-morphisms by sending a bimodule homomorphism to the canonical natural transformation it induces.
A version of this $(\infty,2)$-functor is discussed in \cite[\textsection 7.1]{FHLT} and \cite[Rmk~3.3]{TelemanLectures}.
\end{ex}

\subsubsection*{$(\infty,n)$-functors of $(\infty,n)$-categories of manifolds}

Functors out of the indexing shapes that we have previously discussed are the underlying structures for several mathematical objects used to formalize ideas
in mathematical physics. Indeed, given a manifold $M$ and (symmetric monoidal) $(\infty,n)$-category $\mathscr{C}$, we have that:
\begin{enumerate}[leftmargin=*]
\item A \emph{factorization algebra} valued in $\mathscr{C}$ in the sense of \cite[\textsection6]{CostelloGwilliam1}
defines in particular an $(\infty,1)$-functor
\[\mathcal{F}\colon\mathscr{O}pen_M\to\mathscr{C}.\]
\item A \emph{topological quantum field theory} valued in $\mathscr{C}$ in the sense of \cite{AtiyahTQFT} defines in particular an $(\infty,1)$-functor
\[Z\colon\mathscr{C}ob_{n-1,n}\to\mathscr{C}.\]
\item An \emph{extended framed topological quantum field theory} valued in $\mathscr{C}$ as in \cite{BaezDolan} or \cite[\textsection 1.2]{luriecobordism} defines in particular an $(\infty,n)$-functor
\[Z\colon\mathscr{C}ob_{0,n,\infty}^{\mathrm{fr}}\to\mathscr{C}.\]
\item An \emph{$\mathbb{E}_n$-algebra} valued in $\mathscr{C}$ as in \cite[\textsection2.7]{AFprimer} (cf.~also \cite[\textsection5.1]{LurieHA}) defines in particular an $(\infty,n)$-functor
\[\mathcal{A}\colon\mathscr{D}isk_n^{\mathrm{fr}}\to\mathscr{C}.\]
\item A \emph{topological chiral homology} valued in $\mathscr{C}$ in the sense of \cite[\textsection5.5.2]{LurieHA} and a \emph{homology theory} in the sense on \cite[\textsection 3.4]{AFprimer}
define in particular an $(\infty,1)$-functor
\[H\colon\mathscr{M}fld^{\mathrm{fr}}_{n}\to\mathscr{C}.\]
\end{enumerate}

\section{The homotopy theory of $(\infty,n)$-categories}
\label{nequi}

\label{Equivalences}

As one would guess, $(\infty,n)$-categories and $(\infty,n)$-functors are the starting point of a homotopy theory for which the objects of interest are the $(\infty,n)$-categories themselves. We discuss how to build the $(\infty,1)$-category of $(\infty,n)$-categories, and its axiomatic description by Barwick--Schommer-Pries \cite{BarwickSchommerPries}.

\subsection{Equivalences between $(\infty,n)$-functors}

Given two $(\infty,n)$-functors $F\colon\mathscr{C} \to\mathscr{D}$ and $G\colon\mathscr{C} \to\mathscr{D}$, an \emph{equivalence} between these $(\infty,n)$-functors $\varphi\colon F\simeq G$ consists of, for all $k\geq0$, a map
\[\varphi_k\colon k\operatorname{Mor}\mathscr{C}\to(k+1)\operatorname{Mor}\mathscr{D}\]
that assign to each $k$-morphism in $\mathscr{C}$ a $(k+1)$-equivalence in $\mathscr{D}$
\[a\quad\mapsto\quad\left[\varphi_k(a)\colon F_ka\simeq G_ka\right].\]
One of the requested conditions is that, given a $k$-morphism $f\colon a\to b$, there exists a $(k+1)$-equivalence
\[Gf\circ\varphi(a)\simeq\varphi(b)\circ Ff.\]
Like in previous situations, one would have to encode higher coherence as part of the structure for such an equivalence. We write $F\simeq G$ if there exists an equivalence between the $(\infty,n)$-functors $F$ and $G$. Each of the models proposed to model the notion of an $(\infty,n)$-category must come with a corresponding notion of an equivalence between $(\infty,n)$-functors that packages the appropriate amount of coherence.

\begin{ex}
Every natural isomorphism between ordinary functors defines an equivalence $(\infty,1)$-functors.
\end{ex}

\begin{ex}
Every homotopy between maps of spaces defines an equivalence of $(\infty,0)$-functors.
\end{ex}

\subsection{The $(\infty,1)$-category of $(\infty,n)$-categories}

The notion of an equivalence of $(\infty,n)$-functors plays the role of a higher morphism in an appropriate higher category of $(\infty,n)$-categories and $(\infty,n)$-functors, which is constructed in \cite[\textsection7,8]{BarwickSchommerPries}.

\begin{thm}
\label{AxiomsBSP}
There is an $(\infty,1)$-category
$(\infty,n)\mathscr{C}at$ that satisfies the axioms from  \cite[\textsection7]{BarwickSchommerPries}.
In the $(\infty,1)$-category $(\infty,n)\mathscr{C}at$, the objects are
the $(\infty,n)$-categories, the morphisms are the $(\infty,n)$-functors, the $2$-morphisms are the equivalences $(\infty,n)$-functors. In particular, for all $(\infty,n)$-categories $\mathscr{C}$ and $\mathscr{D}$, there is an $(\infty,0)$-category $\mathscr{F}un(\mathscr{C},\mathscr{D})$ of $(\infty,n)$-functors from $\mathscr{C}$ to $\mathscr{D}$.
\end{thm}

While the precise statement of the axioms from  \cite[\textsection7]{BarwickSchommerPries} is difficult, we can give an informal idea of some of them.
First, it is expected that the collection of $k$-morphisms in an $(\infty,n)$-category 
can be corepresented by an $(\infty,n)$-category, called \emph{$k$-cell}.
More generally, certain collections of composable $k$-morphisms with specified compositions should be corepresented by more general  $(\infty,n)$-categories, which we may refer to as \emph{generalized $k$-cells}. 
One of the requirements is that any $(\infty,n)$-category can be built out of those $k$-cells, and a further requirement encodes a list of relations between the generalized cells. The last requirement which is easy to formulate informally is the fact that $k$-cells detect equivalences between $(\infty,n)$-categories.

We will see in \cref{complicialthm,nCSSthm,Enrichedthm} that, using \cref{Unicity}, the $(\infty,1)$-category $(\infty,n)\mathscr{C}at$ can be described using different models, too.

\subsection{Equivalences of $(\infty,n)$-categories}

An $(\infty,n)$-functor $F\colon\mathscr{C}\to\mathscr{D}$ is an \emph{$(\infty,n)$-equivalence} -- and we write  $F\colon\mathscr{C}\simeq\mathscr{D}$ -- if there exists an $(\infty,n)$-functor $G\colon\mathscr{D}\to\mathscr{C}$ and equivalences of $(\infty,n)$-functors 
\[\operatorname{id}_\mathscr{C}\simeq G\circ F\quad\text{ and }\quad\operatorname{id}_\mathscr{D}\simeq F\circ G.\] We write $\mathscr{C}\simeq\mathscr{D}$ if there exists an equivalence between $\mathscr{C}$ and $\mathscr{D}$. Each of the models proposed to model the notion of an $(\infty,n)$-category must come with a corresponding notion of an $(\infty,n)$-equivalence between $(\infty,n)$-categories.

\begin{ex}
Every equivalence of categories is an $(\infty,1)$-equivalence, and every biequivalence of bicategories is an $(\infty,2)$-equivalence.
\end{ex}

\begin{ex}
Every homotopy equivalence of spaces is an equivalence of $(\infty,0)$-categories.
\end{ex}

The following criterion, which is discussed in \cite[\textsection5.5,5.6]{GH} and in 
\cite[\textsection2.2.1]{ORsurvey}, highlights the inductive nature with respect to $n$ of the notion of an $(\infty,n)$-equivalence between $(\infty,n)$-categories.

\begin{prop}
Let $n>0$. An $(\infty,n)$-functor $F\colon \mathscr{C}\to \mathscr{D}$ between $(\infty,n)$-categories is an equivalence if and only if the following conditions hold.
\begin{enumerate}[leftmargin=*]
    \item The $(\infty,n)$-functor $F$ is \emph{surjective on objects up to equivalence}, that is, for every object
    $b\in\operatorname{Ob}\mathscr{D}$ there exists an object $a\in\operatorname{Ob}\mathscr{C}$ and an $n$-equivalence
    \[b\simeq Fa\text{ in the $(\infty,n)$-category }\mathscr{D}.\]
    \item The $(\infty,n)$-functor $F$ is a \emph{hom-wise equivalence}, that is, for all objects $a,a'\in\operatorname{Ob} \mathscr{C}$ the morphism $F$ induces equivalences
    \[F_{a,a'}\colon \mathscr{C}(a,a')\simeq\mathscr{D}(Fa,Fa')\text{ of $(\infty,n-1)$-categories}.\]
\end{enumerate}
\end{prop}

\subsection{Two examples}

We have mentioned that there are many models to formalize the heuristic idea of an $(\infty,n)$-category, and we can now elaborate a little bit on how one argues that these models are equivalent. More precisely, each model comes with the construction of an $(\infty,1)$-category of $(\infty,n)$-categories, and the models are understood as equivalent if the corresponding $(\infty,1)$-categories of $(\infty,n)$-categories are equivalent. Some of these equivalences \cite{br1,br2,LurieGoodwillie,LoubatonEqui,DKM,ara} have been produced explicitly in the form of Quillen equivalences between model categories. Others can be obtained abstractly using the following recognition principle, which is proven as \cite[Thm~11.2]{BarwickSchommerPries}.

\begin{thm}
\label{Unicity}
An $(\infty,1)$-category $\mathscr{C}$ satisfies the properties from \cite[\textsection7]{BarwickSchommerPries} if and only if there is an $(\infty,1)$-equivalence
\[\mathscr{C}\simeq(\infty,n)\mathscr{C}at.\]
\end{thm}

The formalism of equivalences of $(\infty,n)$-categories can also be fruitfully used to phrase in a compact and complete way many theorems that express correspondences of various kinds. We'll see many instances of this once we have introduced the language of symmetric monoidal $(\infty,n)$-categories, but we can already mention one now.

The following result, which is established as \cite[Thm~3.17]{NikolausSchreiberStevensonGeneral}, is an enhanced classification of principal $G$ bundles over a fixed space $X$, which takes into account how an isomorphism of $G$-bundles corresponds appropriately to an $1$-equivalence in the classifying space $\mathscr{B}G$.
\begin{thm}[Classification of $G$-bundles]
\label{ClassificationGBundles}
Given a space $X$ and a finite group $G$, there is an equivalence
of $(\infty,0)$-categories
\[\mathscr{B}un^G_X\simeq\mathscr{F}un(\mathscr{X},\mathscr{B}G).\]
\end{thm}

\section{Mathematical implementations of $(\infty,n)$-categories}

\label{MathematicalImplementations}

\label{Models}

We discuss in this section some of the established models to formalize mathematically $(\infty,n)$-categories (alongside the corresponding notion of an $(\infty,n)$-functor and of an equivalence of $(\infty,n)$-functors). These models will be:
\begin{itemize}[leftmargin=*]
    \item $n$-complicial sets \cite{VerityComplicialII,or}
    \item $n$-fold complete Segal spaces \cite{BarwickThesis}
    \item categories enriched over a nice model of $(\infty,n-1)$-categories \cite{br1}.
\end{itemize}
See the first section for a list of references on the models that we are not treating here. While we try to make the content as self-contained as possible, in this section we shall assume a basic familiarity with simplicial sets and ordinary category theory. We encourage the reader to consult the various original sources for a complete and precise treatment.

\subsection{Complicial sets}

We recall the definition of Verity's $n$-complicial sets \cite[\textsection 1.3]{or}, based on \cite{VerityComplicialI}, and briefly describe what the axioms are supposed to encode. We refer the reader to \cite{EmilyNotes} for a more gentle explanation of the intuition behind this model. This approach is an enhancement of Joyal's quasi-categories \cite{JoyalVolumeII} and Lurie's naturally marked quasi-categories
\cite[\textsection3.1]{htt}  for $n=1$ and of Lurie's $\infty$-bicategories \cite[\textsection4.2]{LurieGoodwillie} for $n=2$.

\subsubsection*{The definition}

\begin{defn}
An \emph{$n$-complicial set} consists of a simplicial set $X$ with a specified set of simplices (which always include the degenerate ones) in positive dimensions, called \emph{marked}, such that the following hold:
\begin{enumerate}[leftmargin=*]
       \item \emph{Horn fillers}: Given the $k$-horn of an $m$-simplex\footnote{The \emph{$k$-horn of an $m$-simplex} consists of the datum of the boundary of a standard $m$-simplex in which furthermore the $k$-th face is removed.} valued in $X$, if all faces containing the vertices $\{k-1,k,k+1\}\cap[m]$ are marked, then there is a marked $m$-simplex filling the horn.
    \item \emph{Thinness extension}:
    Given a marked $m$-simplex in $X$ in which all faces containing the set of vertices $\{k-1,k,k+1\}\cap[m]$ are marked, if the $(k-1)$st face and the $(k+1)$st faces are marked, then the $k$-th face is also marked.
    \item \emph{Saturation:} Given a marked $m$-simplex, if all its faces containing $\{0,2\}$ and all faces containing $\{1,3\}$ are marked, then all its faces are marked.
    \item \emph{Triviality:} Every simplex in dimension higher than $n$ is marked.
\end{enumerate}
\end{defn}

The notion of an $n$-complicial set is designed to implement the notion of an $(\infty,n)$-category, based on the following interpretation.
\begin{center}
\begin{tabular}{ p{3.5cm} c p{5cm} }
$n$-complicial set & $\leftrightsquigarrow$ & $(\infty,n)$-category \\ 
\hline\\[-0.35cm]
 $0$-simplex & $\leftrightsquigarrow$ & object \\  
 $1$-simplex & $\leftrightsquigarrow$ & morphism \\
 marked $1$-simplex&$\leftrightsquigarrow$&$1$-equivalence\\
 $2$-simplex &$\leftrightsquigarrow$& $2$-morphism with a specified target factorization\\
 marked $2$-simplex &$\leftrightsquigarrow$&composable pair of morphisms with specified composite\\
 $k$-simplex&$\leftrightsquigarrow$&$k$-morphism with certain boundary decomposition\\
 marked $k$-simplex&$\leftrightsquigarrow$&$k$-equivalence with certain boundary decomposition\\
 degenerate $k$-simplex&$\leftrightsquigarrow$&identity morphism with certain boundary decomposition
\end{tabular}
\end{center}

The axioms are designed, respectively, to guarantee the following:
\begin{enumerate}[leftmargin=*]
\setcounter{enumi}{-1}
    \item Every identity is an equivalence.
    \item Given the datum of some composable $k$-morphisms, there exists a composite $k$-morphism.
    \item The composite of $k$-equivalences is a $k$-equivalence.
    \item Every $k$-equivalence admits an inverse.
    \item Every $k$-morphism for $k>n$ is a $k$-equivalence.
\end{enumerate}

\subsubsection*{The $(\infty,1)$-category of $n$-complicial sets}

Denote by $\Delta[1]_t$ the standard $1$-simplex endowed with the maximal marking.

\begin{defn}
A \emph{map of $n$-complicial sets} is a simplicial map preserving the marking.
A map $F\colon X\to Y$ of $n$-complicial sets is a \emph{weak equivalence of $n$-complicial sets} if there exists $G\colon Y\to X$ as well as marking preserving simplicial maps between marked simplicial sets
\[H\colon \Delta[1]_t\times X\to X\text{ and }K\colon \Delta[1]_t\times Y\to Y\]
such that
\[H|_{\{0\}\times X}=G\circ F,\quad H|_{\{1\}\times X}=\operatorname{id}_X\quad K|_{\{0\}\times Y}=F\circ G\quad K|_{\{1\}\times Y}=\operatorname{id}_Y\]
\end{defn}

In this context, an \emph{$(\infty,n)$-category} is an $n$-complicial set, an \emph{$(\infty,n)$-functor} is then just a map of $n$-complicial sets and an \emph{equivalence of $(\infty,n)$-categories} is a weak equivalence of $n$-complicial sets.

The following is a consequence of \cref{nCSSthm} together with the comparison constructed in \cite{LoubatonEqui}.

\begin{thm}
\label{complicialthm}
There is an $(\infty,1)$-category $n\mathscr{C}omp$ in which the objects are $n$-complicial sets and the $1$-morphisms are the maps of $n$-complicial sets,
and which satisfies the axioms from \cite[\textsection7]{BarwickSchommerPries}. In particular, by \cref{Unicity} there is an equivalence of $(\infty,1)$-categories
\[n\mathscr{C}omp\simeq (\infty,n)\mathscr{C}at.\]
\end{thm}

\begin{rmk}
The advantages for the model of $(\infty,n)$-categories given by $n$-complicial sets include the following.
\begin{itemize}[leftmargin=*]
    \item \emph{Combinatorial nature:} The combinatorial nature of the model makes it similar in various aspects to the lower dimensional cases of Kan complexes and quasi-categories, allowing one to reuse several techniques from simplicial homotopy theory. This model has been by now extensively studied, and several aspects of the theory have already been developed in detail.
    \item \emph{Difficulty doesn't grow with $n$:} The nature of $n$-complicial sets is such that working with large $n$ does not present substantial extra complexity compared to the cases of low values of $n$.
    \item \emph{Easy equivalences of $(\infty,n)$-functors:} The presence of marking allows one to recognize immediately which $k$-simplices represent $k$-equivalences.
    \item \emph{Easy $(\infty,n)$-functors:} The definition of map of $n$-complicial sets is very easy to work with, since it's just a marking preserving map.
\end{itemize}
By contrast, some disadvantages include the following.
\begin{itemize}[leftmargin=*]
    \item \emph{Each $k$-morphism is overly represented:} As soon as $k>1$, there are a priori multiple $k$-simplices representing the same $k$-morphism, one for each boundary decomposition that fits the shape of the simplex boundary. For instance, the $2$-complicial set representing a $2$-category with exactly one non-identity $2$-morphism has infinitely many non-degenerate simplices in its implementation as a $2$-complicial set (see \cite[\textsection 2]{NerveSuspension}).
    \item \emph{Composition of $k$-morphisms is hard to access:} As soon as $k>1$, in order to compose two $k$-morphisms one has to rearrange the datum of the corresponding $k$-simplices in a way that fits an appropriate horn, and then evoke the existence of a lift for this horn.
\end{itemize}
\end{rmk}


\subsection{$n$-Fold complete Segal spaces}

We recall the definition of Barwick's $n$-fold complete Segal spaces \cite{BarwickThesis}, and briefly describe what the axioms are supposed to encode. We refer the reader to \cite{CalaqueScheimbauer,HaugsengSpan} for a more gentle explanation of the intuition. 
This approach recovers the one of Rezk's complete Segal spaces \cite{rezk} for $n=1$, and is related to the approach of Bergner--Rezk's categories weakly enriched in complete Segal objects \cite{br1}.

\subsubsection*{The definition}

The following is consistent with \cite[\textsection14]{BarwickSchommerPries}, and is equivalent to the original definition \cite[\textsection 2]{BarwickThesis}.
Recall from \cite[\textsection 2.3]{BarwickThesis} (under the name of Rezk-$\mathcal{M}$-categories for $\mathcal{M}$ being a (nice) model category)
the notion of a \emph{complete Segal object} in a model category.
\begin{defn}
A \emph{$1$-fold complete Segal space} is a complete Segal space in the sense of \cite[\textsection4,6]{rezk}, namely a complete Segal object $X\colon\Delta^{\operatorname{op}}\to\mathit{s}\mathcal{S}\!\mathit{et}_{(\infty,0)}$ in the model structure $\mathit{s}\mathcal{S}\!\mathit{et}_{(\infty,0)}$ for $(\infty,0)$-categories.
For $n>1$, an \emph{$n$-fold complete Segal space} consists of an $n$-fold simplicial space $X\colon(\Delta^n)^{\operatorname{op}}\to\mathit{s}\mathcal{S}\!\mathit{et}$ such that
\begin{enumerate}[leftmargin=*]
    \item For all $i\geq0$, the functor $X_{i,\bullet,\dots,\bullet}\colon(\Delta^{n-1})^{\operatorname{op}}\to\mathit{s}\mathcal{S}\!\mathit{et}_{(\infty,0)}$ is an $(n-1)$-fold complete Segal space.
    \item The functor $X\colon \Delta^{\operatorname{op}} \to \mathit{s}\mathcal{S}\!\mathit{et}^{(\Delta^{n-1})^{\operatorname{op}}}_{(\infty,n-1)}$
    is a complete Segal object in the model structure $\mathit{s}\mathcal{S}\!\mathit{et}^{(\Delta^{n-1})^{\operatorname{op}}}_{(\infty,n-1)}$ for $(\infty,n-1)$-categories from \cite[\textsection2.3]{BarwickThesis}.
    \item The functor $X_{0,\bullet,\dots,\bullet}\colon(\Delta^{n-1})^{\operatorname{op}}\to\mathit{s}\mathcal{S}\!\mathit{et}_{(\infty,0)}$ is homotopically discrete $(\infty, n-1)$-category, in the sense of \cite[\textsection 14]{BarwickSchommerPries}.
\end{enumerate}
\end{defn}

Taking the viewpoint that an \emph{$(\infty,n)$-category} is an $n$-fold complete Segal space, the axioms are designed to support the following intuition (cf.~\cite[Remark 3.13]{HaugsengMorita})
\begin{center}
\begin{tabular}{ p{1.7cm} c p{7cm} }
$n$-fold CSS & $\leftrightsquigarrow$ & $(\infty,n)$-category \\ 
\hline\\[-0.35cm]
 $X_{0,\dots,0}$ & $\leftrightsquigarrow$ & space of objects \\  
$X_{1,\dots,0}$ & $\leftrightsquigarrow$ & space of $1$-morphisms \\
 $X_{1,1,0\dots,0}$ &$\leftrightsquigarrow$& space of $2$-morphisms\\
 $X_{1,\dots,1,0\dots,0}$ &$\leftrightsquigarrow$& space of $k$-morphisms\\
 $X_{1,\dots,1}$ &$\leftrightsquigarrow$&space of $n$-morphisms\\
 $X_{2,0,\dots,0}$&$\leftrightsquigarrow$&space of pairs of composable $1$-morphisms with a specified composite\\
$X_{k,0,\dots,0}$&$\leftrightsquigarrow$&space of $k$-tuples of composable $1$-morphisms with a specified composite\\
\end{tabular}
\end{center}

\subsubsection*{The $(\infty,1)$-category of $n$-fold complete Segal spaces}

\begin{defn}
A \emph{map of $n$-fold complete Segal spaces} is a multisimplicial map.
A map of $n$-fold complete Segal spaces is a \emph{weak equivalence of $n$-fold complete Segal spaces} if it is a levelwise equivalence of spaces.
\end{defn}

In this context, an \emph{$(\infty,n)$-category} is $n$-fold complete Segal space, an \emph{$(\infty,n)$-functor} is then just a map of $n$-fold complete Segal spaces and an \emph{equivalence of $(\infty,n)$-categories} is a weak equivalence of $n$-fold complete Segal spaces.

The following is a consequence of \cite[\textsection 14]{BarwickSchommerPries}.

\begin{thm}
\label{nCSSthm}
There is an $(\infty,1)$-category $n\mathscr{C}SS$ in which the objects are $n$-fold complete Segal spaces and the $1$-morphisms are the maps of $n$-fold complete Segal spaces,
and which satisfies the axioms from \cite[\textsection7]{BarwickSchommerPries}. In particular, by \cref{Unicity} there is an equivalence of $(\infty,1)$-categories
\[n\mathscr{C}SS\simeq (\infty,n)\mathscr{C}at.\]
\end{thm}

\begin{rmk}
The advantages for the model of $(\infty,n)$-categories given by $n$-fold complete Segal spaces include the following.
\begin{itemize}[leftmargin=*]
    \item \emph{Good to implement examples with non-concrete morphisms:} Most $(\infty,n)$-ca\-te\-go\-ries of cobordisms, iterated spans, or flavors of higher Morita categories are typically implemented as $n$-fold complete Segal spaces.
    \item \emph{Prone to model independence:} The $(\infty,1)$-category of $n$-fold complete Segal spaces is understood as a localization of the $(\infty,1)$-category of space-valued presheaves over $\Delta^n$. As such, this can be nicely implemented using a specific model of spaces (most notably Kan complexes), or without reference to a specific model of spaces, working model independently.
    \item \emph{Easy $(\infty,n)$-functors and $(\infty,n)$-equivalences:} The definition of map and weak equivalence of $n$-fold complete Segal spaces are very easy to work with, since it's just a map of presheaves, resp.~ a levelwise weak equivalence.
\end{itemize}
By contrast, some disadvantages include the following.
\begin{itemize}[leftmargin=*]
    \item \emph{Completeness is extremely subtle:} While it is often possible to produce explicit examples of iterated Segal spaces, the completeness condition is often hard to achieve \emph{by hand}. This often forces one to abstractly complete a given presheaf to obtain an $n$-fold complete Segal space, and the process of completion, which exists for abstract reasons is rather uncontrolled (for instance, it generally affects every level of the multisimplicial object). This features for instance in the formalization of the various $(\infty,n)$-categories of cobordisms \cite{luriecobordism,CalaqueScheimbauer}.
\item \emph{$k$-Equivalences hard to access:} It is not clear how to recognize whether a $k$-morphism is a $k$-equivalence.
\end{itemize}
\end{rmk}


\subsection{Categories enriched over $(\infty,n-1)$-categories}

Assuming we have a nice model of $(\infty,n-1)$-categories, there is a model of $(\infty,n)$-categories based on categories strictly enriched over $(\infty,n-1)$-categories; see \cite[\textsection A.3]{htt} and \cite[\textsection 3.10]{br1}. This approach generalizes Bergner's simplicial categories \cite{bergner} and to some extent 
Lack's $2$-categories \cite{lack1,lack2}.

\subsubsection*{The definition}

Suppose to have a good notion of $(\infty,n-1)$-category, $(\infty,n-1)$-functor and equivalence of $(\infty,n-1)$-categories with a well-behaved cartesian product (in technical terms, an excellent cartesian closed model structure for $(\infty,n-1)$-categories). For $n-1=0$, one could take the model structure for Kan complexes as a nice model for $(\infty,0)$-categories. For $n-1=1$, one could take the model structure for quasi-categories or the one for complete Segal spaces as a nice model for $(\infty,1)$-categories. For $n-1=2$, one could take $2$-complicial sets as a nice model for $(\infty,2)$-categories. For general $n-1$ one could take $(n-1)$-complicial sets or complete Segal $\Theta_{n-1}$-spaces from \cite{rezkTheta} as a nice model for $(\infty,n-1)$-categories. For the remainder of the subsection we will assume that one such nice model for $(\infty,n-1)$-categories (with corresponding notions of $(\infty,n-1)$-functor, equivalence of $(\infty,n-1)$-functors,
etc.) has been fixed.

\begin{defn}
A \emph{category enriched over $(\infty, n-1)$-categories} $\mathcal{C}$ consists of
\begin{itemize}[leftmargin=*]
    \item a set of objects
    \[\operatorname{Ob}\mathcal{C};\]
    \item for all $c,c'\in\operatorname{Ob}\mathcal{C}$, an $(\infty,n-1)$-category \[\mathcal{C}(c,c');\]
    \item for all $c,c',c''\in\operatorname{Ob}\mathcal{C}$ 
an $(\infty,n-1)$-functor
    \[\circ\colon\mathcal{C}(c,c')\times \mathcal{C}(c',c'')\to\mathcal{C}(c,c'');\]
    \item for all $c\in\operatorname{Ob}\mathcal{C}$ an $(\infty,n-1)$-functor
    \[\operatorname{id}_c\colon\top\to\mathcal{C}(c,c),\]
    where $\top$ denotes the terminal $(\infty,n-1)$-category.
\end{itemize}
It is required that the following axioms be met.
\begin{itemize}[leftmargin=*]
    \item \emph{Associativity:} for all $c,c',c''\in\operatorname{Ob}\mathcal{C}$ the diagram of $(\infty,n-1)$-categories commutes
    \[
\begin{tikzcd}[column sep=2cm]
\mathcal{C}(c,c')\times\mathcal{C}(c',c'')\times\mathcal{C}(c'',c''')\arrow[r,"{\operatorname{id}_{\mathcal{C}(c,c')}\times\circ}"]\arrow[d, "{\circ\times\operatorname{id}_{\mathcal{C}(c'',c''')}}" swap]&\mathcal{C}(c,c')\times\mathcal{C}(c',c''')\arrow[d,"\circ"]\\
\mathcal{C}(c,c'')\times\mathcal{C}(c'',c''')\arrow[r,"\circ" swap]&\mathcal{C}(c,c''')
\end{tikzcd}
\]
    \item \emph{Unitality:} for all $c,c'\in\operatorname{Ob}\mathcal{C}$ the diagrams of $(\infty,n-1)$-categories commute
\[
\begin{tikzcd}
\top\times\mathcal{C}(c,c')\arrow[dr, "\operatorname{pr}_2"] \arrow[d, "{\operatorname{id}_c \times \operatorname{id}_{\mathcal{C}(c,c')}}"swap]&\\
\mathcal{C}(c,c)\times\mathcal{C}(c,c')\arrow[r,"\circ" swap]&\mathcal{C}(c,c')
\end{tikzcd}
\mbox{ and }\quad
\begin{tikzcd}
\mathcal{C}(c,c')\times \top\arrow[dr, "\operatorname{pr}_1"]\arrow[d, "{\operatorname{id}_{\mathcal{C}(c,c')}\times \operatorname{id}_{c'}}" swap]&\\
\mathcal{C}(c,c')\times\mathcal{C}(c',c')\arrow[r,"\circ" swap]&\mathcal{C}(c,c')
\end{tikzcd}
\]
\end{itemize}
\end{defn}

The model is designed so that the following intuition holds.
\begin{center}
\begin{tabular}{ >{\raggedleft\arraybackslash}p{3cm} c p{4cm} }
cat enriched over $(\infty,n-1)$-cats&$\leftrightsquigarrow$&$(\infty,n)$-category\\
    \hline\\[-0.35cm]
     object& $\leftrightsquigarrow$ & object\\
   object of hom&$\leftrightsquigarrow$& $1$-morphism\\
    morphism of hom&$\leftrightsquigarrow$&$2$-morphism\\
    $k$-morphism in hom&$\leftrightsquigarrow$& $(k+1)$-morphism\\
\end{tabular}
\end{center}

\subsubsection*{The $(\infty,1)$-category of categories enriched over $(\infty,n-1)$-categories}

\begin{defn}
\label{EnrichedFunctor}
A \emph{functor enriched over $(\infty,n-1)$-categories} consists of
\begin{itemize}[leftmargin=*]
    \item map on objects
    \[F\colon\operatorname{Ob}\mathcal{C}\to\operatorname{Ob}\mathcal{D}\]
    \item an $(\infty,n-1)$-functor on hom-$(\infty,n-1)$-categories
    \[F_{c,c'}\colon\mathcal{C}(c,c')\to\mathcal{D}(Fc,Fc')\]
\end{itemize}
such that
\begin{itemize}[leftmargin=*]
    \item \emph{Composition:} $F$ is compatible with composition, that is, for every $c,c',c''\in\operatorname{Ob}\mathcal{C}$ of the diagram of $(\infty,n-1)$-categories commutes 
\[
\begin{tikzcd}[column sep=2cm]
\mathcal{C}(c,c')\times \mathcal{C}(c',c'')\arrow[r,"F_{c,c'}\times F_{c',c''}"]\arrow[d, "\circ"']&\mathcal{D}(Fc,Fc')\times \mathcal{D}(Fc',Fc'')\arrow[d, "\circ"]\\
\mathcal{C}(c,c'')\arrow[r,"F_{c,c''}" swap]&\mathcal{D}(Fc,Fc'')
\end{tikzcd}
\]
    \item \emph{Identity:} $F$ is compatible with identity, that is, for every $c\in\operatorname{Ob}\mathcal{C}$ of the diagram of $(\infty,n-1)$-categories commutes
\[
\begin{tikzcd}
&\top\arrow[ld,"{\operatorname{id}_c}" swap]\arrow[rd,"\operatorname{id}_{Fc}"]&\\
\mathcal{C}(c,c)\arrow[rr,"F_{c,c}" swap]&&\mathcal{D}(Fc,Fc)
\end{tikzcd}
\]
\end{itemize}
\end{defn}

\begin{defn}
A \emph{weak equivalence of categories enriched over $(\infty,n-1)$-categories} is a functor $F\colon\mathcal{C}\to\mathcal{D}$ enriched over $(\infty,n-1)$-categories such that
\begin{itemize}[leftmargin=*]
\item the $(\infty,n)$-functor $F$ is \emph{surjective on objects up to equivalence}, that is, for every object
    $b\in\operatorname{Ob}\mathcal{D}$ there exists an object $a\in\operatorname{Ob}\mathcal{C}$ and an $n$-equivalence
    \[b\simeq Fa\text{ in the $(\infty,n)$-category }\mathcal{D}.\]
    \item the $(\infty,n)$-functor $F$ is a \emph{hom-wise equivalence}, that is, for all objects $a,a'\in\operatorname{Ob} \mathcal{C}$ the morphism $F$ induces equivalences
    \[F_{a,a'}\colon \mathcal{C}(a,a')\simeq\mathcal{D}(Fa,Fa')\text{ of $(\infty,n-1)$-categories}.\]
\end{itemize}
\end{defn}

The following is a consequence of \cref{nCSSthm} together with the comparisons constructed in \cite{br1,br2}.

\begin{thm}
\label{Enrichedthm}
There is an $(\infty,1)$-category $\mathscr{C}at_{(\infty,n-1)\mathscr{C}at}$ in which the objects are the categories enriched over $(\infty, n-1)$-categories, the $1$-morphisms include
the functors enriched over $(\infty, n-1)$-categories,
and which satisfies the axioms from \cite[\textsection7]{BarwickSchommerPries}. In particular, by \cref{Unicity} there is an equivalence of $(\infty,1)$-categories
\[\mathscr{C}at_{(\infty,n-1)\mathscr{C}at}\simeq (\infty,n)\mathscr{C}at.\]
\end{thm}

\begin{rmk}
The advantages for the model of $(\infty,n)$-categories given by categories strictly enriched over a nice model of $(\infty,n-1)$-categories include the following.
\begin{itemize}[leftmargin=*]
    \item \emph{Rather intuitive:} Assuming one has a good understanding of the model of $(\infty,n-1)$-categories that is being used, working with categories strictly enriched in that model reflects the most intuitive take one would have to approach $(\infty,n)$-categories.
    \item \emph{Morphisms are nicely stored:} The inductive nature of the definition allows one to match somewhat exactly the cells of each hom-$(\infty,n-1)$-category with cells of the $(\infty,n)$-category.
    \item \emph{Top dimensional composition is strict:} It is helpful in some contexts that the top dimensional composition is strictly defined, strictly associative and strictly unital.
    \item \emph{Good for concrete examples:} The strictly enriched models are good to implement $(\infty,n)$-categories of higher categories of some sort, often coming from model categories.
\end{itemize}
By contrast, some disadvantages include the following.
\begin{itemize}[leftmargin=*]
    \item \emph{Bad for non-concrete examples:} It is almost impossible to produce an explicit implementation of examples involving cobordisms, spans or higher Morita categories as a category strictly enriched over $(\infty,n-1)$-categories.
    \item \emph{The definition of an $(\infty,n)$-functor is very complicated:} The simple and naive definition of $(\infty,n)$-functor from \cref{EnrichedFunctor} is too restrictive, and insufficient to capture all $(\infty,n)$-functors. To fix this one would have to first take a cofibrant replacement of the domain, a notably painful process.
\end{itemize}
\end{rmk}

\section{Symmetric monoidal $(\infty,n)$-categories}
\label{monoidal}

In this section we introduce the reader to the notion of a symmetric monoidal $(\infty,n)$-category from a heuristic viewpoint, and discuss a list of examples. We point the reader to external sources for a detailed mathematical treatment.

\subsection{Schematic definition}

We start by describing the type of structure expected from a symmetric monoidal $(\infty,n)$-category, and afterwards build towards the extra (structural) axioms that need to be added.

\subsubsection*{The heuristic structure of a symmetric monoidal $(\infty,n)$-category}
\label{AxiomsSymMon}
A \emph{symmetric monoidal $(\infty,n)$-category} $\mathscr{C}$ consists of an $(\infty,n)$-category $\mathscr{C}$ together with a functor
\[\otimes\colon\mathscr{C}\times\mathscr{C}\to\mathscr{C},\quad (f,g)\mapsto f\otimes g.\]
and an object $I_{\mathscr{C}}$ in $\mathscr{C}$.
In particular, this entails that, given any $k$-morphisms $f\colon a\to a'$ and $f'\colon b\to b'$, there is a $k$-morphism $f\otimes f'\colon a\otimes a'\to b\otimes b'$
for all $k\geq0$.

The first layer of structural axioms that one should require for the symmetric monoidal $(\infty,n)$-category $\mathscr{C}$ would be the following:
\begin{itemize}[leftmargin=*]
\item
\emph{Associativity:} Given $k$-morphisms $f,g,h$
    there is a $(k+1)$-equivalence
\[f\otimes (g\otimes h)\simeq(f\otimes g)\otimes h.\]
\item \emph{Unitality:} Given a $k$-morphism
$f$, there are $(k+1)$-equivalences
\[f\otimes\operatorname{id}_{I_\mathscr{C}}\simeq f\quad\text{ and }\quad \operatorname{id}_{I_\mathscr{C}}\otimes f\simeq f.\]
    \item \emph{Symmetry:} Given $k$-morphisms 
    $f,g$,
    there is a $(k+1)$-equivalence
    \[f\otimes g\simeq g\otimes f.\]
\end{itemize}

 Like in previous situations, one would have to further record an infinite amount of coherence data.
 There are different approaches to capturing the notion of a symmetric monoidal $(\infty,n)$-category in a precise mathematical framework. These include e.g.\ \cite[\textsection3]{BarwickThesis}, \cite[\textsection 3]{CalaqueScheimbauer}, \cite[Def.\ 2.0.0.7]{LurieHA}, \cite[\textsection 6.2]{GH}.

\subsection{Examples}

We discuss several examples of symmetric monoidal $(\infty,n)$-ca\-te\-go\-ries that appear naturally for various values of $n$.

\subsubsection*{Structural examples of symmetric monoidal $(\infty,n)$-categories}

We collect here examples of symmetric monoidal $(\infty,n)$-categories that arise from other known variants of (higher) categories.

\begin{ex}
Every symmetric monoidal category can be naturally regarded as a symmetric monoidal $(\infty,1)$-category. This is formalized in
\cite[\textsection2]{LurieHA}.
\end{ex}
\begin{ex}
Every symmetric monoidal bicategory
in the sense of \cite[\textsection3.1]{SPthesis} or \cite{OsornoSymMonBicat} (see also \cite[\textsection3.2]{SPthesis} for more history details) can be naturally regarded as a symmetric monoidal $(\infty,2)$-category. This could be formalized using e.g.~\cite[App.~S]{NikolausScholze}.
\end{ex}

\begin{ex}
Given $n\geq0$, every symmetric monoidal $(\infty,n)$-category can be naturally regarded as a symmetric monoidal $(\infty,n+1)$-category with the same objects, and it can also be regarded as an $(\infty,n+1)$-category with exactly one object by appropriately shifting the dimension of the morphisms. This is formalized in \cite[\textsection6.3]{GH}.
\end{ex}

\subsubsection*{Symmetric monoidal $(\infty,n)$-categories of structured sets}

We give a family of examples of symmetric monoidal $(\infty,n)$-categories in which the objects are structured sets and the morphisms are appropriate structure-preserving functions.

\begin{ex}[Spaces]
The $(\infty,1)$-category $\mathscr{S}pc$ can be enhanced to a symmetric monoidal $(\infty,1)$-category $\mathscr{S}pc^\times$ via the cartesian product
\[\times\colon\mathscr{S}pc\times\mathscr{S}pc\to\mathscr{S}pc.\]
The monoidal unit is the singleton.
The symmetric monoidal $(\infty,1)$-category $\mathscr{S}pc^\times$ can be understood as the underlying $(\infty,1)$-category of the cartesian model structure from \cite[Prop.\ 4.2.8]{hovey}.
\end{ex}

\begin{ex}[Vector spaces and Hilbert spaces]
Given a field $\mathbb{K}$, the $(\infty,1)$-category $\mathscr{V}ect_{\mathbb{K}}$ can be enhanced to a symmetric monoidal $(\infty,1)$-category -- in fact symmetric monoidal category -- $\mathscr{V}ect_{\mathbb{K}}^{\otimes_\mathbb{K}}$ via the ordinary $\mathbb{K}$-tensor product
\[\otimes_\mathbb{K}\colon\mathscr{V}ect_\mathbb{K}\times\mathscr{V}ect_\mathbb{K}\to\mathscr{V}ect_\mathbb{K}.\]
The monoidal unit is the ground field $\mathbb{K}$. Similarly, the $(\infty,1)$-category $\mathscr{H}ilb_{\mathbb{K}}$ can be enhanced to a symmetric monoidal $(\infty,1)$-category $\mathscr{H}ilb_{\mathbb{K}}^{\otimes_\mathbb{K}}$ via the $\mathbb{K}$-tensor product of Hilbert spaces from \cite[\textsection 3.4]{WeidmannHilbert}.
\end{ex}

\begin{ex}[Chain complexes]
Given a field $\mathbb{K}$, the $(\infty,1)$-category $\mathscr{C}h_{\mathbb{K}}$ can be enhanced to a symmetric monoidal $(\infty,1)$-category $\mathscr{C}h_{\mathbb{K}}^{\widehat\otimes_\mathbb{K}}$ via the derived $\mathbb{K}$-tensor product
\[\widehat\otimes_\mathbb{K}\colon\mathscr{C}h_\mathbb{K}\times\mathscr{C}h_\mathbb{K}\to\mathscr{C}h_\mathbb{K}.\]
The monoidal unit is the ground field $\mathbb{K}$ concentrated in degree 0. The symmetric monoidal $(\infty,1)$-category $\mathscr{C}h_{\mathbb{K}}^{\widehat\otimes_\mathbb{K}}$ can be understood as the underlying $(\infty,1)$-category of the symmetric monoidal model structure from \cite[Prop.\ 4.2.13]{hovey} (see also~\cite[Rmk 7.1.2.12]{LurieHA}).
\end{ex}

\begin{ex}[Linear categories]
Given a field $\mathbb{K}$, the $(\infty,2)$-category $\mathscr{C}at^L_{\mathbb{K}}$ can be enhanced to a symmetric monoidal $(\infty,2)$-category $(\mathscr{C}at^L_{\mathbb{K}})^{\boxtimes_\mathbb{K}}$ via an appropriate enhancement of the Deligne--Kelly $\mathbb{K}$-tensor product
\[\boxtimes_\mathbb{K}\colon\mathscr{C}at^L_\mathbb{K}\times\mathscr{C}at^L_\mathbb{K}\to\mathscr{C}at^L_\mathbb{K}\]
from \cite[\textsection 3.2]{BenZviBrochierJordanIntegrating}.
The monoidal unit is the linear category of vector spaces.
\end{ex}

\subsubsection*{Symmetric monoidal $(\infty,n)$-categories with unusual morphisms}

We give a family of examples of symmetric monoidal $(\infty,n)$-categories in which the objects are \emph{not} structured sets and the morphisms are \emph{not} structure-preserving functions.

\begin{ex}[Classifying space of abelian group]
Given an abelian group $G$ with operation $*$, the $(\infty,0)$-category $\mathscr{B}G$ can be enhanced to a symmetric monoidal $(\infty,0)$-category $\mathscr{B}G^*$ via the $(\infty,0)$-functor
\[*\colon\mathscr{B}G\times\mathscr{B}G\to\mathscr{B}G\]
induced by the operation $*$.
\end{ex}

\begin{ex}[Spaces and spans]
Given $n\geq0$, an appropriate version of the $(\infty,n)$-category $\mathscr{S}pan_n$ can be enhanced to a symmetric monoidal $(\infty,n)$-category $(\mathscr{S}pan_n)^\times$ via the cartesian product
\[\times\colon\mathscr{S}pan_n\times\mathscr{S}pan_n\to\mathscr{S}pan_n.\]
The monoidal unit is the singleton. The symmetric monoidal $(\infty,n)$-category $(\mathscr{S}pan_n)^\times$ is discussed in \cite[\textsection1.2, Thm~1.2]{HaugsengSpan}.
\end{ex}

\begin{ex}[Algebras and Bimodules]
Given a field $\mathbb{K}$, the $(\infty,2)$-category $\mathscr{M}or_{\mathbb{K}}$ can be
enhanced to a symmetric monoidal $(\infty,2)$-category $\mathscr{M}or_{\mathbb{K}}^{\otimes_\mathbb{K}}$ via the ordinary $\mathbb{K}$-tensor product
\[\otimes_\mathbb{K}\colon\mathscr{M}or_\mathbb{K}\times\mathscr{M}or_\mathbb{K}\to\mathscr{M}or_\mathbb{K}.\]
The monoidal unit is the ground field $\mathbb{K}$. The symmetric monoidal $(\infty,2)$-category $\mathscr{M}or_{\mathbb{K}}^{\otimes_\mathbb{K}}$ is discussed in \cite[\textsection 4.6]{SPthesis}.
\end{ex}

\subsubsection*{Symmetric monoidal $(\infty,n)$-categories of manifolds}

We mention the traditional enhancement to symmetric monoidal $(\infty,n)$-categories for the $(\infty,n)$-categories of manifolds that we discussed previously.

\begin{ex}[Manifolds, disks]
Given $n\geq0$, the $(\infty,1)$-category $\mathscr{M}fld^{\mathrm{fr}}_{n}$ can be enhanced to a symmetric monoidal $(\infty,1)$-category $(\mathscr{M}fld^{\mathrm{fr}}_{n})^\amalg$ via the disjoint union
\[\amalg\colon\mathscr{M}fld^{\mathrm{fr}}_{n}\times\mathscr{M}fld^{\mathrm{fr}}_{n}\to\mathscr{M}fld^{\mathrm{fr}}_{n}.\]
The monoidal unit is the empty manifold. Similarly, the $(\infty,1)$-category $(\mathscr{D}isk_n^{\mathrm{fr}})^{\amalg}$
can be enhanced to a symmetric monoidal $(\infty,1)$-category $\mathscr{D}isk_n^{\amalg}$. 
The symmetric monoidal $(\infty,1)$-categories $(\mathscr{M}fld^{\mathrm{fr}}_{n})^\amalg$ and $(\mathscr{D}isk_n^{\mathrm{fr}})^{\amalg}$ are discussed in \cite[\textsection2.4,2.6]{AFprimer}.
\end{ex}

\begin{ex}[Open subsets]
Given a space $M$, the $(\infty,1)$-category $\mathscr{O}pen_M$ can be enhanced to a symmetric monoidal $(\infty,1)$-category $\mathscr{O}pen_M^\cup$ via the union
\[\cup\colon\mathscr{O}pen_M\times\mathscr{O}pen_M\to\mathscr{O}pen_M.\]
The monoidal unit is the empty set.
\end{ex}

\begin{ex}[Cobordisms]
Given $n\geq0$, the various versions of the $(\infty,1)$-category $\mathscr{C}ob_{n-1,n}$ can be enhanced to a symmetric monoidal $(\infty,1)$-category -- in fact symmetric monoidal category $(\mathscr{C}ob_{n-1,n})^\amalg$ via the disjoint union
\[\amalg\colon\mathscr{C}ob_{n-1,n}\times\mathscr{C}ob_{n-1,n}\to\mathscr{C}ob_{n-1,n},\]
and similarly the various versions of the $(\infty,n)$-category $\mathscr{C}ob_{0,n,\infty}$ can be enhanced to a symmetric monoidal $(\infty,n)$-category $(\mathscr{C}ob_{0,n,\infty})^\amalg$. The monoidal unit is the empty manifold.
For instance, the symmetric monoidal category $(\mathscr{C}ob_{n-1,n}^{\mathrm{or}})^\amalg$ is mentioned in \cite[Ex.\ 3.2.44]{Kock} and the symmetric monoidal $(\infty,n)$-category $\mathscr{C}ob_{0,n,\infty}^{\amalg}$ is discussed in \cite[\textsection 7]{CalaqueScheimbauer}.
\end{ex}

\section{Symmetric monoidal $(\infty,n)$-functors}
\label{symmonfun}

In this section we introduce the reader to the notion of a symmetric monoidal $(\infty,n)$-functor from a heuristic viewpoint, and discuss a list of examples. We point the reader to external sources for a detailed mathematical treatment.

\subsection{Schematic definition}

We start by describing the type of structure expected from a symmetric monoidal $(\infty,n)$-functor, and afterwards build towards the extra (structural) axioms that need to be added.

\subsubsection*{The heuristic structure}
A \emph{symmetric monoidal $(\infty,n)$-functor}
\[F\colon\mathscr{C}^{\otimes}\to\mathscr{D}^{\boxtimes}\]
consists of an $(\infty,n)$-functor
$F\colon\mathscr{C}\to\mathscr{D}$
satisfying appropriate structural axioms.
The first instance is the requirement that for all $k$-morphisms $f$ and $g$ there exists a $(k+1)$-equivalence in $\mathscr{D}$
\[F(f\otimes g)\simeq Ff\boxtimes Fg.\]
One would have to then add further higher coherence requirements, expressing in particular the compatibility of this $(k+1)$-equivalence with the ones from \cref{AxiomsSymMon}.

\subsection{Examples}

We discuss several examples of $(\infty,n)$-functors that appear naturally for various values of $n$.

\subsubsection*{Explicit examples}

Here are some examples of $(\infty,n)$-functors involving some of the explicit examples of symmetric monoidal $(\infty,n)$-categories seen so far.

\begin{ex}
\label{ExampleBunFunctor2}
The $(\infty,2)$-functor from \cref{ExampleBunFunctor} should
 define a symmetric monoidal $(\infty,2)$-functor
\[\mathscr{B}un^G_{(-)}\colon(\mathscr{C}ob_{0,2,\infty}^{\mathrm{or}})^\amalg\to(\mathscr{S}pan_2)^\times,\]
based on the canonical identification of $(\infty,0)$-categories, a.k.a., spaces 
\[\mathscr{B}un^G_{X\amalg X'}\simeq \mathscr{B}un^G_X\times \mathscr{B}un^G_{X'}.\]
\end{ex}

\begin{ex}
\label{ExamplePathAlgebra2}
Given a field $\mathbb{K}$, for appropriate versions of the symmetric monoidal $(\infty,2)$-categories $(\mathscr{S}pan_2)^\times$ and $\mathscr{M}or_{\mathbb{K}}^{\otimes_\mathbb{K}}$, there should be a symmetric monoidal $(\infty,2)$-functor
\[(\mathscr{S}pan_2)^\times\to\mathscr{M}or_{\mathbb{K}}^{\otimes_\mathbb{K}}\]
that assigns to each space its path algebra (see \cite[\textsection3]{TelemanLectures} and the case $n=2$ of \cite[\textsection 8]{FHLT}).
\end{ex}

\begin{ex}
\label{ExampleCategorification2}
Given a field $\mathbb{K}$, for appropriate versions of the symmetric monoidal $(\infty,2)$-categories $\mathscr{M}or_{\mathbb{K}}^{\otimes_\mathbb{K}}$ and $\mathscr{C}at_{\mathbb{K}}^{\boxtimes_{\mathbb{K}}}$, the $(\infty,2)$-functor from \cref{ExampleCategorification} should define a symmetric monoidal $(\infty,2)$-functor
\[\mathcal{M} od\colon\mathscr{M}or_{\mathbb{K}}^{\otimes_{\mathbb{K}}}\to\mathscr{C}at_{\mathbb{K}}^{\boxtimes_{\mathbb{K}}},\]
based on the canonical identification of linear categories 
\[\mathcal{M} od_{A\otimes_\mathbb{K}A'}\simeq\mathcal{M} od_A\boxtimes_\mathbb{K}\mathcal{M} od_{A'}.\]
\end{ex}

\subsubsection*{Examples from other contexts}

Functors out of the indexing shapes that we have previously discussed are the underlying structure for several mathematical objects used to formalize ideas from mathematical physics. Indeed, given an appropriate manifold $M$ and a symmetric monoidal $(\infty,n)$-category $\mathscr{C}^{\otimes}$, we have that:
\begin{enumerate}[leftmargin=*]
    \item A symmetric monoidal $(\infty,n)$-functor 
    \[\mathcal{F}\colon\mathscr{O}pen_M^\cup\to\mathscr{C}^{\otimes}\]
satisfying a further cosheaf-like condition is the underlying object of a \emph{factorization algebra} valued in $\mathscr{C}^{\otimes}$ in the sense of \cite[\textsection6]{CostelloGwilliam1}.
\item A symmetric monoidal functor
\[Z\colon\mathscr{C}ob_{n-1,n}^\amalg\to\mathscr{C}^{\otimes}\] is precisely a \emph{topological quantum field theory} valued
in a symmetric monoidal category $\mathscr{C}^{\otimes}$ in the sense of \cite{AtiyahTQFT}.
\item A symmetric monoidal $(\infty,n)$-functor
\[Z\colon(\mathscr{C}ob_{0,n,\infty}^{\mathrm{fr}})^\amalg\to\mathscr{C}^{\otimes}\]
is precisely the underlying structure for an \emph{extended framed topological quantum field theory} valued in
$\mathscr{C}^{\otimes}$ as in \cite[\textsection 1.2]{luriecobordism}.
\item A symmetric monoidal $(\infty,n)$-functor
\[\mathcal{A}\colon(\mathscr{D}isk_n^{\mathrm{fr}})^\amalg\to\mathscr{C}^{\otimes}\]
is precisely an \emph{$\mathbb{E}_n$-algebra} valued in $\mathscr{C}^{\otimes}$ in the sense of \cite[\textsection2.7]{AFprimer} (cf.~also \cite[\textsection5.1.1]{LurieHA}).
\item A symmetric monoidal $(\infty,n)$-functor
\[H\colon(\mathscr{M}fld^{\mathrm{fr}}_{n})^\amalg\to\mathscr{C}^{\otimes}\]
that satisfies an appropriate excision property is precisely a \emph{homology theory} valued in $\mathscr{C}^{\otimes}$ in the sense of \cite[\textsection3.4]{AFprimer} and a \emph{topological chiral homology} in the sense of \cite[\textsection5.5.2]{LurieHA}.
\end{enumerate}

Let's mention some of the motivating examples for some of these structures:

\begin{ex}
If $X$ is a pointed space, its $n$-fold loop space defines an $(\infty,n)$-functor, in fact an $\mathbb{E}_n$-algebra valued in $\mathscr{S}pc^\times$,
\[\Omega^nX\colon(\mathscr{D}isk_n^{\mathrm{fr}})^\amalg\to\mathscr{S}pc^\times,\]
as proven in \cite[\textsection 5]{MayIterated}.
\end{ex}

\begin{ex}
Given an $\mathbb{E}_n$-algebra $\mathcal{A}\colon\mathscr{D}isk_n^{\mathrm{fr}}\to\mathscr{C}^{\otimes}$, \emph{factorization homology} of 
$\mathcal{A}$ defines a symmetric monoidal $(\infty,1)$-functor
\[\int_{-}\mathcal{A}\colon(\mathscr{M}fld^{\mathrm{fr}}_{n})^\amalg\to\mathscr{C}^\otimes,\]
as discussed in \cite[\textsection 3.1]{AFprimer} and \cite[\textsection5.5.2]{LurieHA}.
\end{ex}

\begin{ex}
Given a finite group $G$, assuming to have at one's disposal appropriate implementations of the symmetric monoidal $(\infty,2)$-functors from \cref{ExampleBunFunctor2,ExamplePathAlgebra2,ExampleCategorification2}, one can construct a fully extended $2$-dimensional topological quantum field theory valued in the symmetric monoidal $(\infty,2)$-category $\mathscr{C}at_{\mathbb{K}}^{\boxtimes_{\mathbb{K}}}$ (or $\mathscr{M}or^{\otimes}$ or $(\mathscr{S}pan_2)^\times$)
\[Z_G\colon(\mathscr{C}ob_{0,2,\infty}^{\mathrm{or}})^\amalg\to(\mathscr{S}pan_2)^\times\to\mathscr{M}or_{\mathbb{K}}^{\otimes_{\mathbb{K}}}\to\mathscr{C}at_{\mathbb{K}}^{\boxtimes_{\mathbb{K}}}\]
by composing the symmetric monoidal $(\infty,2)$-functors from these previous examples. 
This topological quantum field theory should model a $2$-dimensional gauge theory with finite gauge group $G$
(cf.~\cite[(3.5)]{TelemanLectures} and \cite[\textsection8]{FHLT}).
\end{ex}

\section{Applications}
\label{symmonequi}

We recall the notion of an equivalence between $(\infty,n)$-functors, which is used to define an appropriate $(\infty,1)$-category of symmetric monoidal functors between symmetric monoidal $(\infty,n)$-categories.
We then discuss several theorems that are appropriately stated in terms of $(\infty,1)$-categories of this kind.

\subsection{Schematic definition}

A symmetric monoidal $(\infty,n)$-functor $F\colon\mathscr{C}\to\mathscr{D}$ is a \emph{symmetric monoidal $(\infty,n)$-equivalence} if it is an equivalence of underlying  $(\infty,n)$-categories (cf.~\cite[Lem.~3.2.2.6]{LurieHA}).

The following $(\infty,0)$-category can be obtained as a mapping space in the $(\infty,1)$-category of commutative algebras in an appropriate $(\infty,1)$-category $(\infty,n)\mathscr{C}at$, as in \cite[\textsection2.1.3]{LurieHA}.

\begin{prop}
Given symmetric monoidal $(\infty,n)$-categories $\mathscr{C}^{\otimes}$ and $\mathscr{D}^{\boxtimes}$, there is an
$(\infty,0)$-category
$\mathscr{F}un^{\otimes}(\mathscr{C}^{\otimes},\mathscr{D}^{\boxtimes})$
in which the objects are the 
symmetric monoidal 
$(\infty,n)$-functors from $\mathscr{C}^{\otimes}$ to $\mathscr{D}^{\boxtimes}$.
\end{prop}

The collection of commutative algebras in a given $(\infty,1)$-category, so in particular in the $(\infty,1)$-category of $(\infty,n)$-categories, can be built itself as an $(\infty,1)$-category, using \cite[\textsection 2.1.3]{LurieHA}. As in any $(\infty,1)$-category, we can consider mapping spaces between the objects, e.g.\ as described in \cite[\textsection 1.2.2]{htt}, and those are the spaces of symmetric monoidal $(\infty,n)$-functors which we would like to consider from now on.

The expectation is that, actually, $\mathscr{F}un^{\otimes}(\mathscr{C}^{\otimes},\mathscr{D}^{\boxtimes})$ should have the structure of an interesting $(\infty,n)$-category, rather than an $(\infty,0)$-category, but we won't be treating it in this generality here.

\subsection{Examples}

\label{finalexamples}

We collect several theorems expressed by means of an equivalence of higher categories.

\subsubsection*{The cobordism hypothesis}

\label{FinalApplications}

Extended topological quantum field theories are phrased appropriately with the language of higher categories, and we recall various classification theorems for them.

\begin{defn}
Given $n>0$ and a symmetric monoidal category $\mathscr{C}^{\otimes}$, the \emph{$(\infty,0)$-category of $n$-dimensional TQFTs} valued in $\mathscr{C}^{\otimes}$ is
\[\mathscr{T}qft_{n-1,n}^{\mathrm{fr}}(\mathscr{C}^{\otimes})\coloneq \mathscr{F}un^{\otimes}((\mathscr{C}ob_{n-1,n}^{\mathrm{fr}})^{\amalg},\mathscr{C}^{\otimes}).\]
The \emph{$(\infty,0)$-category $\mathscr{T}qft_{n-1,n}^{\mathrm{or}}(\mathscr{C}^{\otimes})$ of oriented $n$-dimensional TQFTs} is defined analogously, replacing $(\mathscr{C}ob_{n-1,n}^{\mathrm{fr}})^{\amalg}$
with $(\mathscr{C}ob_{n-1,n}^{\mathrm{or}})^{\amalg}$.
\end{defn}

The following is folklore (see e.g.~\cite[Thm~3.2]{CarquevilleRunkel}). 

\begin{prop}[$1$-dimensional cobordism hypothesis]
Given a field $\mathbb{K}$, there is an equivalence of $(\infty,0)$-categories 
\[\mathscr{T}qft^{\mathrm{or}}_{0,1}(\mathscr{V}ect_\mathbb{K}^{\otimes_\mathbb{K}})\simeq \mathscr{V}ect^{\mathrm{fg}}_{\mathbb{K}},\]
where $\mathscr{V}ect^{\mathrm{fg}}_{\mathbb{K}}$ denotes 
the $(\infty,0)$-category of finite-dimensional vector spaces and isomorphisms.
\end{prop}

The following is folklore (see e.g.~\cite[\textsection3.3]{Kock}).

\begin{prop}[$2$-dimensional cobordism hypothesis]
Given a field $\mathbb{K}$, there is an equivalence of $(\infty,0)$-categories
\[\mathscr{T}qft_{1,2}^{\mathrm{or}}(\mathscr{V}ect_\mathbb{K}^{\otimes_\mathbb{K}})\simeq \mathscr{A}lg^{\mathrm{Frob,com}}_{\mathbb{K}}\]
where $\mathscr{A}lg^{\mathrm{Frob,com}}_{\mathbb{K}}$ denotes the $(\infty,0)$-category of commutative Frobenius $\mathbb{K}$-algebras, described e.g.~in \cite[\textsection2]{Kock}, and isomorphisms.
\end{prop}

\subsubsection*{Extended cobordism hypothesis}

We recall various versions of cobordism hypotheses that classify certain fully extended topological quantum field theories that are phrased appropriately with the language of higher categories.

\begin{defn}
Given a symmetric monoidal $(\infty,n)$-category $\mathscr{C}^{\otimes}$ and $n>0$, the \emph{$(\infty,0)$-category of fully extended $n$-dimensional framed TQFTs} valued in $\mathscr{C}^{\otimes}$ is
\[\mathscr{T}qft_{0,n}^{\mathrm{ext,fr}}(\mathscr{C}^{\otimes})\coloneq \mathscr{F}un^{\otimes}((\mathscr{C}ob_{0,n,\infty}^{\mathrm{fr}})^{\amalg},\mathscr{C}^{\otimes}).\]
\end{defn}

The following -- in fact a stronger version of it -- is proven in \cite[\textsection 4.6]{SPthesis}.

\begin{thm}[Extended $2$-dimensional oriented cobordism hypothesis]
Given a field $\mathbb{K}$, there is an equivalence of $(\infty,0)$-categories
\[\mathscr{T}qft_{0,2}^{\mathrm{or}}(\mathscr{M}or_{\mathbb{K}}^{\otimes_{\mathbb{K}}})\simeq \mathscr{A}lg^{\mathrm{Frob,sep}}_{\mathbb{K}}.\]
Here, $\mathscr{A}lg^{\mathrm{Frob,sep}}_{\mathbb{K}}$ denotes an $(\infty,0)$-category of separable Frobenius $\mathbb{K}$-algebras, defined in \cite[\textsection A.2]{SPthesis}, obtained from the symmetric monoidal bicategory from \cite[\textsection4.6]{SPthesis}.
\end{thm}

The following was originally formulated in \cite{BaezDolan}, sketched in \cite{luriecobordism}, and further developed in \cite{AFcobordism,CalaqueScheimbauer,GradyPavlov}, and other ongoing projects.

\begin{thm}[Extended $n$-dimensional framed cobordism hypothesis]
For $n>0$ there is an equivalence of $(\infty,0)$-categories
\[\mathscr{T}qft^{\mathrm{ext,fr}}_{0,n}(\mathscr{C}^{\otimes})\simeq (\mathscr{C}^{\otimes})^{\simeq},\]
where $(\mathscr{C}^{\otimes})^{\simeq}$ denotes the symmetric monoidal sub $(\infty,0)$-category of fully dualizable objects in $\mathscr{C}^{\otimes}$, described in \cite[\textsection 2.3]{luriecobordism}.
\end{thm}

On the physics side, there appears to be a connection between gapped phases of matter and certain kinds of (extended) topological quantum field theories which is rather new; a first insight can be obtained in \cite{GaiottoJF}.

\subsubsection*{Classifications of $\mathbb{E}_n$-algebras}

We recall several classifications theorems of various $(\infty,0)$-categories
of $\mathbb{E}_n$-algebras that are phrased appropriately with the language of higher categories.

\begin{defn}[{\cite[\textsection5.1]{LurieHA}}]
Given a symmetric monoidal $(\infty,1)$-category $\mathscr{C}^{\otimes}$ and $n>0$, the \emph{$(\infty,0)$-category of $\mathbb{E}_n$-algebras} valued in $\mathscr{C}^{\otimes}$ is
\[\mathscr{A}lg_{\mathbb{E}_n}(\mathscr{C}^{\otimes})\coloneq \mathscr{F}un^{\otimes}((\mathscr{D}isk_n^{\mathrm{fr}})^{\amalg},\mathscr{C}^{\otimes}).\]
\end{defn}

The following is folklore (see e.g.~\cite[Ex.~5.1.0.7,Ex.~5.1.1.7]{LurieHA}).

\begin{prop}
Given field $\mathbb{K}$, there is an equivalence of $(\infty,0)$-categories
\[\mathscr{A}lg_{\mathbb{E}_1}(\mathscr{V}ect_\mathbb{K}^{\otimes_\mathbb{K}})\simeq\mathscr{A}lg_\mathbb{K},\]
and for $n>1$ there is an equivalence of $(\infty,0)$-categories
\[\mathscr{A}lg_{\mathbb{E}_n}(\mathscr{V}ect_\mathbb{K}^{\otimes_\mathbb{K}})\simeq\mathscr{A}lg_\mathbb{K}^{\mathrm{com}}\]
where $\mathscr{A}lg^{\mathrm{com}}_{\mathbb{K}}$ denotes the $(\infty,0)$-category of commutative $\mathbb{K}$-algebras and isomorphisms.
\end{prop}

The following -- in fact a stronger version of it -- is essentially \cite[Thm~5.4.5.9]{LurieHA}.

\begin{thm}[Classification of locally constant factorization algebras in $\mathbb{R}^n$]
Given $n>0$ and a symmetric monoidal $(\infty,1)$-category $\mathscr{C}^{\otimes}$ there is an equivalence of $(\infty,0)$-categories
\[\mathscr{F}act\mathscr{A}lg_{\mathbb{R}^n}^{\mathrm{l.c.}}(\mathscr{C}^{\otimes})\simeq\mathscr{A}lg_{\mathbb{E}_n}(\mathscr{C}^{\otimes})\]
where $\mathscr{F}act\mathscr{A}lg_{\mathbb{R}^n}^{\mathrm{l.c.}}$ denotes the $(\infty,0)$-category of locally constant factorization algebras over $\mathbb{R}^n$ valued in $\mathscr{C}^{\otimes}$, as obtained from the symmetric monoidal $(\infty,1)$-category from \cite[Thm~5.4.5.9]{LurieHA}.
\end{thm}

The following -- in fact a stronger version of it -- is proven in \cite[\textsection3.6]{AFprimer}.

\begin{thm}[Classification of excisive homology theories]
Given $n>0$ and a symmetric monoidal $(\infty,1)$-category $\mathscr{C}^{\otimes}$ that is $\otimes$-presentable in the sense of \cite[Def.\ 3.1]{AFprimer}, there is an equivalence of $(\infty,0)$-categories
\[\mathscr{A}lg_{\mathbb{E}_n}(\mathscr{C}^{\otimes})\simeq\mathscr{F}un^{\otimes,\mathrm{exc}}((\mathscr{M}fld^{\mathrm{fr}}_{n})^\amalg,\mathscr{C}^{\otimes})\]
where $\mathscr{F}un^{\otimes,\mathrm{exc}}((\mathscr{M}fld^{\mathrm{fr}}_{n})^\amalg,\mathscr{C}^{\otimes})$ denotes an $(\infty,0)$-category of symmetric monoidal functors $(\mathscr{M}fld^{\mathrm{fr}}_{n})^\amalg\to\mathscr{C}^{\otimes}$ that are excisive in the sense of \cite[Def.~3.28]{AFprimer}.
\end{thm}

\section*{Summary and further reading}
We gave a brief introduction to $(\infty,n)$-categories, functors between them, the correct notion of `sameness' between them, as well as symmetric monoidal analogs of all these notions. We illustrated with examples motivated from mathematical physics how the use of $(\infty,n)$-categories and their symmetric monoidal variants has become ubiquitous in this area, e.g.\ when studying questions related to the topological quantum field theories and gauge theories. There are many further applications of this formalism that we did not treat in this chapter, and which we encourage the reader to explore.

Given that higher category theory is a rather new field and might prove to be challenging to learn from scratch, we end by mentioning a few introductory references that the interested reader could consult to start learning more about some aspects of the mathematical treatment of the subject:

\begin{itemize}[leftmargin=*]
\item \cite{GrothCourse}: an introduction to quasi-categories, which are a model for $(\infty,1)$-categories;
\item \cite{RasekhCSS}: an introduction to complete Segal spaces, which are a model for $(\infty,1)$-categories;
\item \cite{BergnerSurvey1}: a survey on some of the models for $(\infty,1)$-categories and how they compare to each other;
\item \cite{RVytm}: an introductory lecture series on a model-independent treatment for the theory of $(\infty,1)$-categories;
\item \cite{EmilyNotes}: an introductory lecture series to $n$-complicial sets, which are a model for $(\infty,n)$-categories;
\item \cite{luriecobordism}: an introduction to $n$-fold complete Segal spaces, which are a model for $(\infty,n)$-categories and how they can be used to formalize the extended cobordism hypothesis;
\item \cite{BergnerModelsn}: a survey on some of the models for $(\infty,n)$-categories and how they compare to each other;
\item \cite{FreedCobSurvey}: a survey on the cobordism hypothesis.
\item \cite{ORsurvey}: a survey on the notion of equivalence between and inside strict and weak $(\infty,n)$-categories.
\end{itemize}

amsalpha
\bibliographystyle{abbrv}
\bibliography{ref2}%

\end{document}

%% file: PreambleChapterProofs.tex
\usepackage{amssymb}
   \usepackage{mathrsfs}
\usepackage{enumitem} 
 \usepackage{array}
 \usepackage{mathtools}

\usepackage{tikz}
\usetikzlibrary{cd}

 \usepackage{hyperref}
\usepackage[capitalise]{cleveref}

 \theoremstyle{plain}   
\newtheorem{thm}{Theorem}[subsection] 
\makeatletter\let\c@thm\c@thm\makeatother

\makeatletter\let\c@cor\c@thm\makeatother

\makeatletter\let\c@lem\c@thm\makeatother
\newtheorem{prop}{Proposition}[subsection]
\makeatletter\let\c@prop\c@thm\makeatother

\makeatletter\let\c@claim\c@thm\makeatother

\makeatletter\let\c@conjecture\c@thm\makeatother

\theoremstyle{definition}

\newtheorem{defn}{Definition}[subsection]
\makeatletter\let\c@defn\c@thm\makeatother

\makeatletter\let\c@const\c@thm\makeatother

\makeatletter\let\c@notn\c@thm\makeatother

\makeatletter\let\c@convention\c@thm\makeatother

\makeatletter\let\c@convention\c@thm\makeatother

\makeatletter\let\c@convention\c@thm\makeatother

\theoremstyle{remark}

\newtheorem{rmk}{Remark}[subsection]
\makeatletter\let\c@rmk\c@thm\makeatother
\newtheorem{ex}{Example}[subsection]
\makeatletter\let\c@ex\c@thm\makeatother

\makeatletter\let\c@observation\c@thm\makeatother

\makeatletter\let\c@warning\c@thm\makeatother

\makeatletter\let\c@digression\c@thm\makeatother

\makeatletter\let\c@answ\c@thm\makeatother

\makeatletter
\let\c@equation\c@thm
\numberwithin{equation}{subsection}
\makeatother

\crefname{lem}{Lemma}{Lemmas}
\crefname{thm}{Theorem}{Theorems}
\crefname{defn}{Definition}{Definitions}
\crefname{notn}{Notation}{Notations}
\crefname{const}{Construction}{Constructions}
\crefname{prop}{Proposition}{Propositions}
\crefname{rmk}{Remark}{Remarks}
\crefname{cor}{Corollary}{Corollaries}
\crefname{equation}{Display}{Displays}
\crefname{ex}{Example}{Examples}
\crefname{thmalph}{Theorem}{Theorems}
\crefname{answ}{Answer}{Answers}
\crefname{question}{Question}{Questions}

\AtBeginDocument{%
   \def\MR#1{}
}